\documentclass[trsc,nonblindrev]{informs3} % current default for manuscript submission

\OneAndAHalfSpacedXI % current default line spacing
\usepackage{algorithm}
\usepackage{algorithmicx}
\usepackage[noend]{algpseudocode}
\usepackage[flushleft]{threeparttable}
\usepackage{wasysym}
\usepackage{multirow}
\usepackage{hyperref}

\usepackage{natbib}
 \bibpunct[, ]{(}{)}{,}{a}{}{,}%
 \def\bibfont{\small}%
 \usepackage{rotating,booktabs}
\TheoremsNumberedThrough     % Preferred (Theorem 1, Lemma 1, Theorem 2)
\EquationsNumberedThrough    % Default: (1), (2), ...
\begin{document}
%%%%%%%%%%%%%%%%

%\setcounter{page}{1}
%\DOI{10.1287/xxxx.0000.0000}%

 \RUNAUTHOR{Lin, Ghaddar, Nathwani}

\RUNTITLE{Electric Vehicle Routing under Time-Variant Electricity Prices}

% \TITLE{Bundling Information Goods of Decreasing Value}
\TITLE{Electric Vehicle Routing with Charging/Discharging under Time-Variant Electricity Prices}

\ARTICLEAUTHORS{%
\AUTHOR{Bo Lin}
\AFF{University of Toronto, 5 King's College Rd, Toronto, ON M5S 3G8, Canada\\
\EMAIL{blin@mie.utoronto.ca}
\URL{}}
\AUTHOR{Bissan Ghaddar}
%David G. Burgoyne Faculty Fellow
\AFF{Ivey Business School, University of Western Ontario, 1255 Western Road, London, ON N6G 0N1, Canada\\ 
\EMAIL{bghaddar@uwaterloo.ca}
\URL{}}
\AUTHOR{Jatin Nathwani}
\AFF{University of Waterloo, 200 University Avenue W., Waterloo, ON N2L 3G1, Canada\\ 
\EMAIL{nathwani@uwaterloo.ca}
\URL{}}
} % end of the block

\ABSTRACT{%
The integration of electric vehicles (EVs) with the energy grid has become an important area of research due to the increasing EV penetration in today's transportation systems. Under appropriate management of EV charging and discharging, the grid can currently satisfy the energy requirements of a considerable number of EVs. Furthermore, EVs can help enhance the reliability and stability of the energy grid through ancillary services such as energy storage. This paper proposes the EV routing problem with time windows under time-variant electricity prices (EVRPTW-TP) which optimizes the routing of an EV fleet that are delivering products to customers, jointly with the scheduling of the charging and discharging of the EVs from/to the grid. The proposed model is a multiperiod vehicle routing problem where EVs can stop at charging stations to either recharge their batteries or inject stored energy to the grid. Given the energy costs that vary based on time-of-use, the charging and discharging schedules of the EVs are optimized to benefit from the capability of storing energy by shifting energy demands from peak hours to off-peak hours when the energy price is lower. The vehicles can recover the energy costs and potentially realize profits by injecting energy back to the grid at high price periods. EVRPTW-TP is formulated as an optimization problem. A Lagrangian relaxation approach and a hybrid variable neighborhood search/tabu search heuristic are proposed to obtain high quality lower bounds and feasible solutions, respectively. Numerical experiments on instances from the literature are provided. The proposed heuristic is also evaluated on a case study of an EV fleet providing grocery delivery at the region of Kitchener-Waterloo in Ontario, Canada. Insights on the impacts of energy pricing, service time slots, range reduction in winter as well as fleet size are presented.
}%

%\KEYWORDS{deterministic inventory theory; infinite linear programming duality; 
%  existence of optimal policies; semi-Markov decision process; cyclic schedule}

\KEYWORDS{Electric vehicle routing; energy storage; sustainable last-mile delivery; mixed integer programming; Lagrangian relaxation, metaheuristics.}
\HISTORY{}

\maketitle
%%%%%%%%%%%%%%%%%%%%%%%%%%%%%%%%%%%%%%%%%%%%%%%%%%%%%%%%%%%%%%%%%%%%%%

% Samples of sectioning (and labeling) in TRSC
% NOTE: (1) \section and \subsection do NOT end with a period
%       (2) \subsubsection and lower need end punctuation
%       (3) capitalization is as shown (title style).
%
%\section{Introduction.}\label{intro} %%1.
%\subsection{Duality and the Classical EOQ Problem.}\label{class-EOQ} %% 1.1.
%\subsection{Outline.}\label{outline1} %% 1.2.
%\subsubsection{Cyclic Schedules for the General Deterministic SMDP.}
%  \label{cyclic-schedules} %% 1.2.1
%\section{Problem Description.}\label{problemdescription} %% 2.

% Text of your paper here

%======================================================================
\section{Introduction} \label{chap:intro}
%======================================================================

% The trend of vehicle electrification
Over the recent years, sustainability has become a paramount global concern. 
In the transportation sector, public and private institutions are attempting to increase the penetration of electric vehicles (EVs) due to their ability to mitigate greenhouse gas emission and their direct impact on reducing particulate matter pollution \citep{boulanger2011vehicle, waraich2013plug}. Companies are also investigating the use of new green technologies due the potential brand benefits given the growing demands for green products \citep{kleindorfer2005sustainable, dekker2012operations}. UPS, FedEx, and Walmart are among the leading companies that have deployed EV fleets in their operations \citep{UPSStrategy}. As a consequence, the past decade has seen a rapid expansion of EV adoption \citep{McKinsey2018}. According to \cite{CanadaNewMotor2018}, the number of new EV registrations has increased from 25,163 in 2014 to 69,010 in 2018. 
As projected by Canada Energy Board \citep{FutureEnergy}, EVs will account for over 60\% of the new motor vehicle registrations in Canada by 2040. 
% For instance, in 2019, Transport Canada announced an investment of 300 million Canadian dollars to a new federal zero-emission vehicle (ZEV) purchase incentive program \cite{CanadaZEVIncentives}. 

% The relationship between EV and the grid, V2G model (an intro about the basic idea)
Not surprisingly, the growing penetration of EVs has a significant impact on the power system. 
Previous studies have shown that, without proper management, the EVs will represent a sizable fraction of the total demand for energy \citep{trivino2019joint, dyke2010impact}. As a result, the gap between peak and off-peak demands will increase, and the ramp requirements will affect the stability and reliability of the power network \citep{villar2012impact}. 
With optimized scheduling of the charging, the existing power system infrastructure can accommodate the energy requirements of a considerable number of EVs \citep{razeghi2016impacts,kintner2007impacts,letendre2008plug}. 
In doing so, the need for installing new capacity, which is expensive, time-consuming, and harmful to the environment, can be minimized \citep{villar2012impact}. Furthermore, EVs can be incorporated into the power grid as a reliable and cost-effective distributed power storage \citep{kempton1997electric}. By optimizing the charging and discharging to the grid, an EV fleet connected to the energy grid can assist to level out peaks in the overall electricity consumption and support the utilization of intermittent renewable energy.
% , thus alleviate the grid pressure. 
% In addition, V2G can also support the utilization of intermittent renewable energy and generate profits for the grid and EV owners if certain business models are applied.

% The implementation of V2G
Although the vehicle-to-grid (V2G) connectivity represents an enticing idea, it nonetheless remains in the pilot stages of development \citep{sovacool2018neglected}, and is mainly focused around a centralized architecture where a controller manages the ancillary services of a large group of EVs that are charging and discharging energy from/to the grid \citep{guille2009conceptual, sortomme2010optimal}. Commercial EV fleet owners such as logistics and e-commerce companies are naturally strong candidates for realizing the benefits of EV integration in the energy grid as they aggregate a large number of EVs. 
% The key element for commercial EV fleet operation is the efficient routing and scheduling of charging and discharging cycles. 
% To the best of our knowledge, no previous study is able to jointly optimize the routing and charging/discharging for an EV fleet as a whole in bidirectional V2G context.

% What we did
This paper considers a delivery service system that operates a fleet of EVs that are primarily used to deliver products to customers. The EVs can be charged or discharged at the home depot or at charging stations in the network. When charging, a cost is paid according to the energy price at the time of use. If energy is discharged to the grid, a profit is creditted to the EV. We thus propose the EV routing problem with time windows under time-variant electricity prices (EVRPTW-TP) which optimizes the monetary cost of the EV fleet operation while allowing the charging and discharging of EVs at time-of-use electricity prices. EVRPTW-TP extends the EV routing problem with time windows \citep{schneider2014electric} by incorporating additional operational constraints, allowing the partial charging and discharging of the EVs, and accounting for the time-varying electricity prices.
In order to solve EVRPTW-TP, a hybrid Variable Neighborhood Search/Tabu Search (VNS/TS) heuristic that can generate high quality feasible solutions efficiently is developed. 
A Lagrangian Relaxation that is solved by a cutting plane approach is also proposed to obtain lower bounds. The results are evaluated using a variation of the widely-used vehicle routing instances of \cite{solomon1987algorithms}. Finally, the model is evaluated using a case study of an EV fleet performing grocery delivery in the Kitchener-Waterloo region in Ontario, Canada. Managerial insights are drawn with respect to electricity pricing, time slots design, winter range reduction, and fleet size.

% The contribution of the research
This paper is the first to investigate the joint optimization of routing and charge/discharge scheduling of multiple EVs under time-variant electricity prices. The proposed model provides operational decisions to support commercial EV fleet operators in order to lower the overall energy costs. The proposed model also offers important implications for policy makers and can assist the power system regulators to better understand the impact of EV fleets on energy markets, and to predict and estimate the market reaction to energy price adjustments. The managerial insights extracted can help policy makers in creating more efficient energy pricing schemes to maximize the environmental and economic benefits from the widespread adoption of commercial EV fleets. 

% The organization of the thesis
The rest of this paper is organized as follows. Following this introductory section, Section \ref{chap:LiteratureReview} reviews the related literature. The proposed problem is formulated in Section~\ref{chap::EVRPTWDPSetting}. The Lagrangian relaxation is presented in Section~\ref{chap:LRforEVRPTWDP} and the proposed VNS/TS heuristic  is then discussed in Section~\ref{chap:VNSTSforEVRPTWDP}. Computational results and the case study are presented in Sections~\ref{compresults} and \ref{chap:UseCase}, respectively. Finally, Section \ref{chap:Conclusion} concludes and highlights future research directions.

% objective

% contribution:
% 1) new business model 2) data-driven framework for commercial EV fleet operation 3) user behavior modeling (dearth of research in this area) 4) in terms of algorithm, integration of ML and OR 

%======================================================================
\section{Literature Review} \label{chap:LiteratureReview}
%======================================================================

The vehicle routing problem (VRP) was first proposed by \cite{dantzig1959truck} as a generalization of the well-known travelling salesman problem. In general, given a set of geographically scattered customers each associated with a demand, the VRP seeks to assign customers to vehicles in such a manner that the demand of each customer is satisfied while the total distance travelled by the fleet of vehicles to serve all the customers is minimized. 

Since the introduction of the classical VRP, numerous variants were developed and investigated to account for realistic constraints and objectives. One of the most common variants is the VRP with time windows where visits to individual customers are restricted to fixed time intervals \citep{russell1977effective}. % Another extention is the electric vehicle routing problems which is a special variant of 
Another variation is the green vehicle routing problem introduced by \cite{erdougan2012green} which particularly models alternative fuel vehicles and accounts for the opportunity to extend a vehicle's distance limitation by visiting an en-route station facility to replenish. \cite{schneider2014electric} tailors this framework specifically to EVs and proposes the EV routing problem with time windows (EVRPTW). Instead of using a constant replenishment time as in \cite{erdougan2012green}, EVRPTW assumes a linear energy charging time that is associated with the battery level of the electric vehicle upon arrival to a station. 

{Following \cite{erdougan2012green} and \cite{schneider2014electric}, various studies have investigated the optimization of EV routing and charging/discharging operations which are summarized in Table \ref{lit_review}. In the uni-directional V2G contexts, \cite{felipe2014heuristic} and \cite{keskin2016partial} consider partial charging strategies with multiple types of chargers, each with a different charging speed and static unit cost. \cite{yang2015electric} optimize over the monetary cost of an EV providing pickup and delivery services under time-variant charging prices. The problem is extended by \cite{barco2017optimal} to a multi-vehicle case incorporating the battery degradation cost. Also under time-variant electricity prices, \cite{yu2018autonomous} models a fleet of autonomous EVs that provide customer delivery and renewable energy storage to the energy grid. A quadraticly-constrained mixed integer program is formulated and three objective functions are proposed to either minimize the total driving distance, to maximize the amount of energy charged from renewable sources, or to minimize the amount of time until the vehicles reach the final destinations.}

{In the bi-directional V2G contexts where EVs are allowed to inject energy back to the grid, \cite{tang2017joint} consider a set of EVs travelling from their origins to the corresponding destinations without en-route customers. EVs can detour for en-route charging and discharging to two types of stations, one providing renewable energy at a low price while the other is a normal station with higher charging cost and provide discharging reward. The energy prices do not vary across time, the objective is to minimize the overall monetary cost of the EV fleet. \cite{trivino2019joint} extend the problem by incorporating intermediate stops for EVs and considering time-variant electricity prices as well as battery degradation cost. However, the model is not able to coordinate the schedules for different EVs and does not have a customer delivery component. \cite{abdulaal2017solving} study a similar problem for an EV fleet considering EV congestion at charging/discharging stations. The customer assignments to different EVs are assumed to be given, the problem thus de-generate to a single EV case. Moreover, the routing and charging decisions are made sequentially, which are likely to be sub-optimal. To the best of our knowledge, no previous research has been conducted to jointly optimize the routing and charging/discharging operations of multiple EVs under time-variant electricity prices.}

{The model proposed in this paper fills the gap between considering energy networks and transportation independently. The proposed model extends the work of \cite{schneider2014electric} to include energy discharging to the grid in a multi-period framework that also accounts for the changing energy prices. Thus the proposed model can be seen as a bi-directional V2G system where the fleet of EVs whose primary purpose is to deliver products to customers can also be used to store and redistribute energy from/to the energy grid. The EV routing problem with time windows under time-variant electricity prices is presented next.}

\begin{table}[!ht]
\centering
\caption{Existing literature about optimization of EV routing and charging/discharging operations} \label{lit_review}
\begin{scriptsize}
{
\begin{tabular}{@{}cccccc@{}}
\toprule
                             & Customer Delivery & Charging    & Discharging & Time-Variant Prices             & Multi-EVs   \\ \midrule
\cite{schneider2014electric} & \checkmark        & full        &             &                                 & \checkmark     \\
\cite{felipe2014heuristic}   & \checkmark        & partial     &             &                                 & \checkmark    \\
\cite{yang2015electric}      & \checkmark        & partial     &             &  \checkmark                     &             \\
\cite{keskin2016partial}     & \checkmark        & partial     &             &                                 & \checkmark  \\
\cite{desaulniers2016exact}  & \checkmark        & partial     &             &                                 & \checkmark  \\
\cite{tang2017joint}         &                   & partial     & \checkmark  &                                 & \checkmark   \\
\cite{abdulaal2017solving}   & \checkmark        & partial     & \checkmark  &  \checkmark                     &             \\ 
\cite{barco2017optimal}      & \checkmark        & partial     &             &  \checkmark                     & \checkmark  \\
\cite{yu2018autonomous}      & \checkmark        & partial     &             &                                 & \checkmark   \\
\cite{trivino2019joint}      &                   & partial     & \checkmark  &  \checkmark                     &             \\
EVRPTW-TP                    & \checkmark        & partial     & \checkmark  &  \checkmark                     & \checkmark  \\
%\cite{PELLETIER2019225}      & \checkmark        & full        &             &                                 & \checkmark     \\
\bottomrule
\end{tabular}
}
\end{scriptsize}
\end{table}
%======================================================================
\section{EV Routing under Time-Variant Electricity Prices} \label{chap::EVRPTWDPSetting}
%======================================================================
% 
% %======================================================================
% \section{Problem Description} \label{sec:ProbDesc}
% %======================================================================
\subsection{Framework and Assumptions}
To formulate EVRPTW-TP, a complete directed graph $G(V_{c,s,od}, E)$ is considered where $V_{c,s,od}$ denotes the set of all nodes and $E$ is the set of all edges.
The nodes are partitioned into three distinct categories: customer nodes, station nodes, and a depot. 
Let $V_{c} = \{1, 2, \dots, N\}$ denote the set of $N$ customers where each node $i$ is associated with a demand $q_i$, a service time $s_i$, and a time window $[e_i, l_i]$ during which an EV should arrive to node $i$. 
Let $V_{s} = \{N+1,\dots, N+S\}$ be the set of $S$ en-route stations where an EV can charge or discharge energy. The depot is denoted by two nodes $0$ and $N+S+1$, i.e. $V_{od}=\{0,N+S+1\}$, where node $0$ is the start of the vehicle route and $N+S+1$ is the end of the route. {Each station or depot node $i\in V_{s}\cup V_{od}$ has a time window $[e_i, l_i]$ from the start to the end of the planning horizon.} Each edge $(i, j)$ is associated with a distance $a_{ij}$ and time $t_{ij}$ that denotes the travel distance and time between nodes $i$ and $j$, respectively. 

{An EV's instantaneous power depends on its mass, acceleration, as well as aerodynamic, rolling, and grade resistances \citep{wu2015electric}. The change in EV mass is negligible for certain applications, e.g. small-package delivery, while the resistances could be captured by the driving speed and acceleration assuming related parameters that are given in \citep{schneider2014electric}. Additionally, the power consumption does not directly translate to the battery's state of charge due to the non-linear battery efficiency and other real-world characteristics such as temperature \citep{hannan2017review}. Incorporating these realistic features could more accurately describe the energy consumption and recharging/discharging processes, yet they will introduce computational burdens. For the sake of keeping the model tractable, following \cite{schneider2014electric} and \citet{desaulniers2016exact}, we assume a constant energy consumption rate $g$ and a constant charging speed $\frac{1}{\alpha}$. The energy consumption for edge $(i, j)$ is thus given by $c_{ij} = ga_{ij}$, the time required to recharge the energy consumed by traveling along edge $(i,j)$ is given by $f_{ij} = \alpha c_{ij}$.}

We assume a commercial EV fleet consisting of $K$ homogeneous EVs, each with a load capacity $Q$ and a battery capacity $C$. The time required to fully recharge the battery from $0$ is defined as $B=\alpha C$. 
At the beginning of the planning horizon, all the EVs are at the depot (node $0$) with a full battery. All EVs should return to the depot (node $N+S+1$) before the end of the planning horizon. 
During each visit to a station or to the depot, an EV can either pay to charge its battery or make profits by injecting energy back to the grid from its battery. 
The charging cost and discharging reward vary according to the time period (time-of-use energy pricing). 
We assume that EVs are allowed to perform either charging or discharging during their station visits, but are not allowed to perform both during the same time period.
% We make this assumption to avoid frequent switches between charging and discharging which were shown to have detrimental effects on battery lifespan.

The planning horizon is formed of $|T|$ consecutive discrete periods, each of length $\delta$. {Each time period $t$ refers to a time interval $[\delta(t-1), \ \delta t)$ and is associated with a buying (from the grid) energy price $p^t_b$ and a selling (to the grid) energy price $p^t_s$. We assume that if an EV is to charge or discharge during a time period, then it has to do so for the full time period. Due to the constant linear charging/discharging assumption, a fixed amount of energy $\frac{\delta}{B}C$ would be charged/discharged during each period, we therefore can calculate a charging cost $P_{re}^t = \frac{\delta}{B}C p_b^t$ and a discharging reward $P_{dis}^t = \frac{\delta}{B} Cp_s^t$ for each period.} EVs are recharged back to full battery capacity during the night using at a lower energy price $p_{night}$. For the ease of presentation, we define $P_{night} = \frac{p_{night}}{\alpha}$ as the cost of charging per unit of time during the night. The notations are summarized in Table \ref{NotationTable}.

The problem formulation that is proposed next jointly optimizes the routes and charging/discharging schedule of the $K$ vehicles by minimizing the net electricity cost given that all customer demands are satisfied. 

\begin{table}
\caption{Summary of the Notation}
\label{NotationTable}
\begin{tabular}{@{}cl@{}}
\toprule
                & \multicolumn{1}{c}{Definition}                                                  \\ \midrule
$V_c$           & Set of customer nodes $V_{c} = \{1, 2, \dots, N\}$                                                    \\
$V_s$           & Set of station nodes $V_{s} = \{N+1,\dots, N+S\}$                                                    \\
$V_{od}$        & Set of depot nodes $V_{od}=\{0,N+S+1\}$                                             \\
$V_{c,s}$     & $V_{c,s} = V_c \cup V_s$\\
$V_{c,s,o}$     & $V_{c,s,o} = V_c \cup V_s \cup \{0\}$\\
$V_{c,s,d}$     & $V_{c,s,d} = V_c \cup V_s \cup \{N+S+1\}$\\
$V_{s,d}$     & $V_{s,d} = V_s \cup \{N+S+1\}$\\
$V_{s,o}$      & $V_{s,o} = V_s \cup \{0\}$\\
$V_{s,od}$      & $V_{s,od} = V_s \cup V_{od}$\\
$V_{c,od}$     & $V_{c,od} = V_c \cup V_{od}$\\
$V_{c,s,od}$     & $V_{c,s,od} = V_c \cup V_s \cup V_{od}$\\
$p$ & Depot node $p = N+S+1$\\
$E$             & Set of all edges                                                            \\ 
$T$             & Set of charging/discharging periods                                                            \\ 
$N$             & Number of customers                                                         \\
$S$             & Number of stations                                                          \\
$K$             & Number of EVs                                                               \\
$\delta$        & Length of each charging/discharging period                                              \\
$Q$             & Cargo capacity                                                       \\
$C$             & Battery capacity                                                                \\
$B$             & Amount of time required to fully charge the EV battery from empty\\
$\alpha$        & The reciprocal of the constant charging speed                       \\
$g$             & Energy consumption rate with respect to distance traveled                        \\
$a_{ij}$        & Travel distance of edge $(i, j)$                                                \\
$t_{ij}$        & Travel time of edge $(i, j)$                                                \\
$c_{ij}$        & Energy consumption along edge $(i, j)$                                         \\
$f_{ij}$        & The amount of time required to charge the energy consumed along edge $(i, j)$ \\
$e_i$           & Earliest service start time at node $i$                                       \\
$l_i$           & Latest service start time at node $i$                                         \\
$s_i$           & Required service time at node $i$                                             \\
$q_i$             & The demand at node $i$                      \\
$p_{night}$     & Unit electricity buying (from the grid) price at night ($\cent$/kWh)                          \\
$p_{b}^t$       & Unit electricity buying (from the grid) price during period $t$ ($\cent$/kWh)                           \\
$p_{s}^t$       & Unit electricity selling (to the grid) price during period $t$ ($\cent$/kWh)                            \\
$P_{night}$     & The cost of charging per unit of time at night                           \\
$P_{re}^{t}$    & Cost of charging during period $t$                                \\
$P_{dis}^{t}$   & Reward for discharging during period $t$                           \\
\bottomrule
\end{tabular}
\end{table}

%======================================================================
\subsection{Problem Formulation} \label{sec:MathFormu}
%======================================================================

% \begin{table}
% \centering
% \caption{Summary of the Decision Variables}
% \label{DecisionVar}
% \begin{tabular}{@{}ccl@{}}
% \toprule
% Variable    & Type       & \multicolumn{1}{c}{Definition}                                                \\ \midrule
% $x_{ijk}$   & Binary     & If edge $(i,j)$ is travelled by vehicle $k$ ($=1$) or not ($=0$)              \\
% $y_{ik}$    & Binary     & Whether vehicle $k$ charges ($=1$) or discharges at node $i$ ($=0$)           \\
% $r_{itk}$   & Binary     & If vehicle $k$ charges at node $i$ during period $t$ ($=1$) or not ($=0$)     \\
% $d_{itk}$   & Binary     & If vehicle $k$ discharges at node $i$  during period $t$ ($=1$) or not ($=0$) \\
% $\tau_{ik}$ & Continuous & The arrival time of vehicle $k$ at node $i$                                   \\
% $b_{ik}$    & Continuous & The battery level of vehicle $k$ on arrival at node $i$                                  \\
% $u_{ik}$    & Continuous & The remaining cargo of vehicle $k$ on arrival at node $i$                                \\ \bottomrule
% \end{tabular}
% \end{table}

To formulate EVRPTW-TP, the following binary variables are introduced:
\begin{align*}
 x_{ijk} &= \begin{cases}
            &1 \quad \mbox{if edge $(i,j)$ is traveled by vehicle $k$,}\\
            &0 \quad \mbox{otherwise,}
           \end{cases}\\
%  y_{ik} &= \begin{cases}
%             &1 \quad \mbox{if vehicle $k$ charges its battery at node $i$,}\\
%             &0 \quad \mbox{otherwise,}
%            \end{cases}\\
r_{itk} &= \begin{cases}
            &1 \quad \mbox{if vehicle $k$ charges its battery at node $i$ at time period $t$,}\\
            &0 \quad \mbox{otherwise,}
           \end{cases}\\
d_{itk} &= \begin{cases}
            &1 \quad \mbox{if vehicle $k$ discharges its battery at node $i$ at time period $t$,}\\
            &0 \quad \mbox{otherwise.}
           \end{cases}
\end{align*}
The continuous variables are 
\begin{align*}
 \tau_{ik} &:\mbox{ arrival time of vehicle $k$ at node $i$,                                  } \\
b_{ik}    &:\mbox{ remaining energy (in terms of charging time) in vehicle $k$ upon arrival to node $i$,}\\
u_{ik} &:\mbox{ remaining cargo in vehicle $k$ upon arrival to node $i$.} 
\end{align*}
The decision variables $r_{itk}$ and $d_{itk}$ are only associated with the station and depot nodes, while the remaining variables are associated with all the nodes. EVRPTW-TP is formulated as 

\begin{align}
 \label{obj:EVRPTWDP}
    \min &\quad \sum_{k=1}^{K}\sum_{i\in V_{s,od}}\sum_{t\in T}
    [r_{itk}P_{re}^{t} - d_{itk}P_{dis}^{t}] + 
    \sum_{k=1}^K P_{night}[B - b_{pk} - 
    \sum_{t\in T}\delta(r_{ptk} - d_{ptk})],\\
\label{con:EVRPTWDP_CustOnce}
\textrm{s.t. }& \quad \sum_{k=1}^K\sum_{j\in V_{c,s,d}} x_{ijk} = 1, 
\quad \forall i \in V_c,\\
\label{con:EVRPTWDP_EVBack}
&\sum_{i\in V_{c,s, o}}x_{ipk} = 1, 
\quad \forall k \in \{1,2,\dots, K\},\\
\label{con:EVRPTWDP_balance}
&\sum_{j\in V_{c, s, o}}x_{jik} - \sum_{j\in V_{c,s,d}}x_{ijk} = 0, 
\quad \forall i \in V_{c,s}, k \in \{1,2,\dots, K\},\\
\label{con:EVRPTWDP_TimeFromCust}
&\tau_{ik} + (t_{ij} + s_{i})x_{ijk} - \delta |T|(1 - x_{ijk}) \leq \tau_{jk}, 
\quad \forall i \in V_{c,s,o}, j \in V_{c,s,d}, k \in \{1,2,\dots, K\},\\
\label{con:EVRPTWDP_TimeFromStat}
%&\tau_{ik} + \sum_{t\in T}\delta (r_{itk} + d_{itk}) + t_{ij}x_{ijk}- M(1 - x_{ijk}) \leq \tau_{jk}, \quad \forall i \in V_{s,o}, j \in V_{c,s,d}, k \in \{1,2,\dots, K\},\\
&t\delta (r_{itk} + d_{itk}) + t_{ij}x_{ijk}- \delta |T|(1 - x_{ijk}) \leq \tau_{jk}, \quad \forall i \in V_{s,o}, j \in V_{c,s,d}, t \in T, k \in \{1,2,\dots, K\},\\
\label{con:EVRPTWDP_TimeWin}
&e_i \leq \tau_{ik} \leq l_i, 
\quad \forall i \in V_{c,s,od}, k \in \{1,2,\dots, K\},\\
\label{con:EVRPTWDP_IniBatt}
&b_{0k} = B, 
\quad \forall k \in \{1,2,\dots, K\},\\
\label{con:EVRPTWDP_BattFromCust}
&b_{jk} \leq b_{ik} - f_{ij}x_{ijk} + B(1-x_{ijk}), 
\quad \forall i \in V_c, j \in V_{c,s,d}, k \in \{1,2,\dots, K\},\\
\label{con:EVRPTWDP_BattFromStat}
&b_{jk} \leq b_{ik} + \sum_{t\in T}\delta r_{itk} - \sum_{t\in T}\delta d_{itk} - f_{ij}x_{ijk} + B(1-x_{ijk}), 
\quad \forall i\in V_{s,o}, j \in V_{c,s,d}, k \in \{1,2,\dots, K\},\\
\label{con:EVRPTWDP_PossCharge}
&\sum_{t\in T}\delta r_{itk} \leq B - b_{ik}, 
\quad \forall i \in V_{s, od}, k \in \{1,2,\dots, K\},\\
\label{con:EVRPTWDP_PossDis}
&\sum_{t\in T}\delta d_{itk} \leq b_{ik}, 
\quad \forall i \in V_{s, od}, k \in \{1,2,\dots, K\},\\
\label{con:EVRPTWDP_ReOrDis1}
&r_{itk} + d_{itk} \leq 1, 
\quad \forall i \in V_{s, od}, t\in T, k \in \{1,2,\dots, K\},\\
% \label{con:EVRPTWDP_ReOrDis1}
% &\sum_{t\in T}r_{itk} \leq |T|y_{ik}, 
% \quad \forall i \in V_{s, od}, k \in \{1,2,\dots, K\},\\
% \label{con:EVRPTWDP_ReOrDis2}
% &\sum_{t\in T}d_{itk} \leq |T|(1 - y_{ik}), 
% \quad \forall i \in V_{s, od}, k \in \{1,2,\dots, K\},\\
\label{con:EVRPTWDP_NoDisBeforeArr}
&\tau_{ik} - (t-1)\delta \leq \delta |T|(1-d_{itk}-r_{itk}), 
\quad \forall i\in V_{s,d}, t\in T, k \in \{1,2,\dots, K\},\\
\label{con:EVRPTWDP_CargoFromVet}
&u_{jk} \leq u_{ik} - q_ix_{ijk} + Q(1-x_{ijk}), 
\quad \forall i \in V_{c,s,o}, j \in V_{c,s,d}, k \in \{1,2,\dots, K\},\\
\label{con:EVRPTWDP_IniCargo}
&u_{0k} = Q, 
\quad k \in \{1,2,\dots, K\},\\
\label{con:EVRPTWDP_BattCon}
&0\leq b_{jk} \leq B\sum_{i\in V_{c,s,o}}x_{ijk},
\quad \forall j \in V_s, k \in \{1,2,\dots, K\},\\
\label{con:EVRPTWDP_ConVar}
&u_{ik}, \ \tau_{ik} \geq 0, 
\quad \forall i \in V_{c,s,od}, k \in \{1,2,\dots, K\},\\
\label{con:EVRPTWDP_ConX}
&x_{ijk} \in \{0, 1\}, 
\quad \forall i \in V_{c,s,o}, \forall j \in  V_{c,s,d}, k \in \{1,2,\dots, K\},\\
% \label{con:EVRPTWDP_ConY}
% &y_{ik} \in \{0, 1\}, 
% \quad \forall i \in V_s, k \in \{1,2,\dots, K\},\\
\label{con:EVRPTWDP_ConRD}
&r_{itk}, \quad d_{itk} \in \{0, 1\}, 
\quad \forall i\in V_{c, od}, t \in T, k \in \{1,2,\dots, K\}.
\end{align}

The objective function \eqref{obj:EVRPTWDP} minimizes the net cost, i.e. the total cost of charging the batteries minus the total reward earned from discharging the batteries. The first part of the objective function corresponds to the net cost during the planning horizon, while the second part is the cost of fully recharging back all EVs at night. 
Constraints \eqref{con:EVRPTWDP_CustOnce} ensure that every customer is served by exactly one EV.
Constraints \eqref{con:EVRPTWDP_EVBack} force all the EVs to return to the depot by the end of the planning horizon. 
Constraints \eqref{con:EVRPTWDP_balance} guarantee that no route ends at a customer or a station node. 
% Constraints \eqref{con:EVRPTWDP_TimeFromCust} -- \eqref{con:EVRPTWDP_TimeWin} account for the arrival time of EVs at a given node.
Constraints \eqref{con:EVRPTWDP_TimeFromCust} ensure time feasibility of the edges leaving the customer and the depot nodes, while constraints~\eqref{con:EVRPTWDP_TimeFromStat} deal with the edges originating from charging stations.
% Note that $t\delta$ is the time when the $t^{th}$ period ends. 
% Constraints~\eqref{con:EVRPTWDP_TimeFromStat} assume that an EV can leave the station right after it completes the charging/discharging. 
Constraints~\eqref{con:EVRPTWDP_TimeWin} ensure that the time windows of all the nodes are not violated. 
Constraints~\eqref{con:EVRPTWDP_IniBatt} indicate that every EV is fully charged before leaving the depot. 
Constraints~\eqref{con:EVRPTWDP_BattFromCust}--\eqref{con:EVRPTWDP_BattFromStat} track the battery capacity along the route. 
Constraints~\eqref{con:EVRPTWDP_PossCharge} indicate that an EV battery cannot be recharged to a level that exceeds its capacity, while Constraints \eqref{con:EVRPTWDP_PossDis} state that an EV battery cannot be discharged to a level below 0. 
Constraints \eqref{con:EVRPTWDP_ReOrDis1} indicate that an EV is allowed to discharge or recharge its battery at a station or the depot nodes but is not allowed to discharge and recharge during the same time period. 
Constraints \eqref{con:EVRPTWDP_NoDisBeforeArr} ensure that an EV cannot start charging/discharging at a station/depot before arrival and before the start of the time period. Constraints \eqref{con:EVRPTWDP_CargoFromVet} ensure that the demands along a route are all satisfied, whereas Constraints \eqref{con:EVRPTWDP_IniCargo} state that all EVs have a full cargo at the start of the planning horizon. 
Constraints \eqref{con:EVRPTWDP_BattCon}--\eqref{con:EVRPTWDP_ConRD} indicate the variables types and limits. 
% {The ``big-Ms'' in time-related constraints \eqref{con:EVRPTWDP_TimeFromCust}, \eqref{con:EVRPTWDP_TimeFromStat}, and \eqref{con:EVRPTWDP_NoDisBeforeArr} should be set as the end of the planning horizon, i.e. $\delta |T|$. They are set to EV battery capacity $B$ in energy-related constraints \eqref{con:EVRPTWDP_BattFromCust} and \eqref{con:EVRPTWDP_BattFromStat}, and EV cargo capacity $Q$ in cargo-related constraints \eqref{con:EVRPTWDP_CargoFromVet}.}

%======================================================================
\section{Lagrangian Relexation for the EVRPTW-TP} \label{chap:LRforEVRPTWDP}
%======================================================================

As shown in Section \ref{sec:performance}, solving EVRPTW-TP to optimality is computationally very challenging. Lagrangian relaxation is a well known algorithm that has been used to address many complex optimization problems. Particularly, in the context of vehicle routing, Lagrangian relaxation has been used to address several variants of VRP \citep{fisher1997vehicle, kallehauge2006lagrangian}. 

For the EVRPTW-TP, all the constraints other than Constraints \eqref{con:EVRPTWDP_CustOnce} are associated with a particular vehicle $k$. Given this special structure which is common in vehicle routing problems, Constraints~\eqref{con:EVRPTWDP_CustOnce} are relaxed and the violation is penalized in the objective function using the Lagrangian multipliers $\lambda_i$. The resulting relaxed problem is
\begin{align*}
% \label{LR_LagRelax}
	Z_{LR}(\lambda) = \min \quad
	&\sum_{k=1}^K \sum_{i\in V_{s, od}}\sum_{t\in T}[r_{itk} P_{re}^{t} - d_{itk} P_{dis}^{t}]  \\ 
	& + \sum_{k=1}^K P_{night}[B - b_{pk} - \sum_{t\in T} \delta(r_{ptk} - d_{ptk})] \\
	& + \sum_{i\in V_c}\lambda_i(1-\sum_{k=1}^K \sum_{j\in V_{c,s,d}}x_{ijk}),\\
	\mbox{s.t. }&\ \eqref{con:EVRPTWDP_EVBack}-\eqref{con:EVRPTWDP_ConRD}.
\end{align*}
Since the EV fleet is homogeneous, problem $Z_{LR}(\lambda)$ decomposes into $K$ identical sub-problems 
\begin{align*}
% \label{LR_SubProb}
    Z_{SP}(\lambda) = \min \quad
    &\sum_{i\in V_{s, od}}\sum_{t\in T} [r_{it}P_{re}^{t} - d_{it}P_{dis}^{t}] + P_{night}[B - b_{p} - \sum_{t\in T}\delta (r_{pt} - d_{pt})] \\
    &- \sum_{i\in V_c} \lambda_i\sum_{j\in V_{c,s,d}}x_{ij},\\
    	\mbox{s.t. }&\ \eqref{con:EVRPTWDP_EVBack}-\eqref{con:EVRPTWDP_ConRD}.
\end{align*}
% Note that the vehicle index $k$ is dropped since $Z_{SP}(\lambda)$ is identical for the $K$ subproblems.

The value of the Lagrangian bound $Z_{LR}(\lambda)$ is 
\begin{align*}
 Z_{LR}(\lambda) = K \times Z_{SP}(\lambda) + \sum_{i\in V_c}\lambda_i,
\end{align*}
and the best Lagrangian bound is given by $\max \limits_{\lambda\in \mathbb{R}^{|V_c|}}Z_{LR}(\lambda)$. Given $H$, the set of feasible solutions of the Lagrangian subproblem, the best Lagrangian bound can be found by solving the following problem where $r^h_{it}$, $d^h_{it}$, and $x^h_{ij}$ describe the charging, discharging, and routing schedule associated with solution $h\in H$, respectively,
\begin{align*}
 \max \limits_{\lambda\in \mathbb{R}^{V_c}} \left\{ K \times \min \limits_{h\in H} \left\{ \sum_{i\in V_{s, od}}\sum_{t\in T} [r^h_{it}P_{re}^{t} - d^h_{it}P_{dis}^{t}] + P_{night}[B - b^h_{p} - \sum_{t\in T}\delta (r^h_{pt} - d^h_{pt})]- \sum_{i\in V_c} \lambda_i\sum_{j\in V_{c,s,d}}x^h_{ij} \right\} + \sum_{i\in V_c}\lambda_i,\right\}
\end{align*}
which is equivalent to the Lagrangian master problem
\begin{align*}
 Z_{MP} = \max  \limits_{\lambda\in \mathbb{R}^{V_c}}  & \sum_{i\in V_c}\lambda_i + K \theta, \\
 \mbox{s.t. } & \theta + \sum_{i\in V_c} \lambda_i\sum_{j\in V_{c,s,d}}x^h_{ij} \leq \sum_{i\in V_{s, od}}\sum_{t\in T} [r^h_{it}P_{re}^{t} - d^h_{it}P_{dis}^{t}] + P_{night}[B - b^h_{p} - \sum_{t\in T}\delta (r^h_{pt} - d^h_{pt})], \ \forall h \in H.
\end{align*}

Since the set $H$ is not known beforehand, we implemented the stablized cutting-plane approach developed by \cite{kallehauge2006lagrangian} to obtain $Z_{MP}$ iteratively starting with an empty set $H$. The determination of the set of cutting planes require solution of the subproblem. Given fixed values of the Lagrangian multipliers $\lambda_i$, the Lagrangian subproblem is solved to obtain $Z_{SP}(\lambda)$ and a new feasible solution $h\in H$. The resulting solution generates a cut that is added to the master problem. The relaxed master problem is solved to obtain new values for the Lagrangian multipliers. The solution of the relaxed master problem provides an upper bound on the optimal Lagrangian bound while the optimal solution of $Z_{LR}(\lambda)$ provides a lower bound. The algorithm iterates until the gap between the upper and lower bounds is sufficiently small. Alternatively, the Lagrangian bound can also be obtained by solving the Dantzig-Wolfe reformulation of $Z_{MP}$ using column generation. Both methods require solving $Z_{SP}(\lambda)$ repetitively.

The Lagrangian sub-problem has similar structure as the column generation sub-problem introduced by \cite{desaulniers2016exact} for EVRPTW. However, it is nontrivial to apply the labeling algorithm proposed by \cite{desaulniers2016exact} to solve $Z_{SP}(\lambda)$ because of the difference between their objective functions. Although the EV's energy consumption along the route is linear with respect to the total distance, the cost of a unit distance in $Z_{SP}(\lambda)$ depends on the time of charging which introduces extra complexity to the problem. In addition, the decision about discharging is largely independent of the travelling distance and hence we can not formulate the Lagrangian sub-problem as a elementary shortest path problem with resource constraints. For these reasons, we use a standard mixed integer programming (MIP) solver to solve $Z_{SP}(\lambda)$. % Developing more efficient algorithms to solve the sub-problem presents interesting opportunities for future research.

%======================================================================
\section{VNS/TS Hybrid Heuristic for EVRPTW-TP} \label{chap:VNSTSforEVRPTWDP}
%======================================================================

The Lagrangian relaxation presented in Section \ref{chap:LRforEVRPTWDP} provides lower bounds on the optimal solution of problem \eqref{obj:EVRPTWDP}--\eqref{con:EVRPTWDP_ConRD}. To obtain upper bounds and very importantly to implement in practice, good quality feasible solutions are needed relatively quickly. Following the framework presented in \cite{schneider2014electric} for EVRPTW, this section presents a hybrid variable neighborhood search and tabu search (VNS/TS) meta-heuristic with an annealing mechanism to solve EVRPTW-TP.

Hybrid VNS/TS meta-heuristics have been previously applied successfully for routing problems \citep{melechovsky2005metaheuristic, tarantilis2008hybrid}. {The overall framework of the VNS/TS is shown in Algorithm~\ref{Alg:HeuristicPseudoCode}.}
VNS/TS consists of three main components: (1) an initialization step which identifies an initial solution, (2) a variable neighborhood search (VNS) to diversify the search process from the current solution, and (3) a tabu search (TS) component that runs for each candidate solution of VNS for local intensification. VNS/TS stops when the iterations limit is reached or when no improving solution is identified after a fixed number of consecutive VNS iterations. The details of each component are presented next.

\begin{algorithm}
    \caption{VNS/TS Heuristic For EVRPTW-TP} 
    \label{Alg:HeuristicPseudoCode}
    \begin{algorithmic}[1]
        \State $S = initialization()$
        \State $counter \gets 0$
        \For {$i = 1, 2, \dots, \eta_{vns}$}
            \State $S^{'} = Move2Neighbor(S)$
            \For{$j = 1, 2, \dots, \eta_{tabu}$}
                \State $S^{''} = Tabu(S^{'})$
            \EndFor
            %\State $S^{''} \gets S^{'}$
            \If{$f_{gen}(S) > f_{gen}(S^{''})$}
                \State $counter \gets 0$
            \Else
                \State $counter \gets counter + 1$
            \EndIf
            \If{$counter \geq \eta_{early}$}
                \State Break
            \ElsIf {$Accept(S, S^{''})$}
                $S \gets S^{''}$          
            \EndIf
        \EndFor
    \end{algorithmic}
\end{algorithm}

%======================================================================
\subsection{Initialization}
%======================================================================
In the initialization step, routes are constructed such that each customer with a positive demand is visited exactly once. For that, the well known sweep heuristic is used to obtain the initial feasible solution {without charging and discharging operations} \citep{cordeau2001unified}. Customers are first sorted according to an increasing order of the geometric angle using as a reference a randomly selected point. Then, starting from the customer with the smallest angle, customers are inserted to the active route at the position resulting in the minimal increase of the travel distance of the route. 
Once the battery or cargo constraints of the active route are violated, a new route is initiated until the number of routes used so far is equal to the EV fleet size. Then, all the remaining customers are inserted into the last route.

%======================================================================
\subsection{Generalized Cost Function}
%======================================================================

The VNS/TS meta-heuristic considers infeasible solutions during the search. Similar to \cite{schneider2014electric}, a cost function is used to evaluate the quality of a solution $S$. The generalized cost function $f_{gen}(S)$ is given by
\begin{equation}
\label{func:GeneralizedCost}
    f_{gen}(S) = f_{elec}(S) + \beta_{tw}\Phi_{tw}(S) + \beta_{batt}\Phi_{batt}(S) + \beta_{cargo}\Phi_{cargo}(S)
\end{equation}
where $f_{elec}(S)$ is the net cost of electricity (charging cost minus discharging reward);  
$\Phi_{tw}(S)$, $\Phi_{batt}(S)$, and $\Phi_{cargo}(S)$ are the violations of the time window, battery, and cargo constraints, respectively; and $\beta_{tw}$, $\beta_{batt}$ and $\beta_{cargo}$ are penalty factors corresponding to each violation. The cumulative net cost $f_{elec}(S)$, and the cumulative violations $\Phi_{tw}(S)$, $\Phi_{batt}(S)$, and $\Phi_{cargo}(S)$ are the sum of the net costs of the individual routes $f_{elec}(R)$ and the individual violations for each route $\Phi_{tw}(R)$, $\Phi_{batt}(R)$, and $\Phi_{cargo}(R)$, respectively.

%======================================================================
\subsubsection{Electricity Cost and Violation Evaluation}
%======================================================================

In order to evaluate the electricity cost and the violations, let $r_i$ denote the $i^{th}$ node along route $R$ of length $n$. The following metrics are then defined for every node along a given route
\begin{equation}
\label{var:EarliestArrival}
    T^{E}_{i} = \left\{
    \begin{aligned}
        & 0,\quad i=1\\
        & \max\left\{\min\left(T^{E}_{i-1}, l_{r_{i-1}}\right)
        + s_{r_{i-1}}
        + t_{r_{i-1}r_i}, %[r_{i-1}, r_{i}], 
        e_{r_i}\right\}, \quad\forall i = 2, 3, \dots, n \\
    \end{aligned}
    \right.
\end{equation}

\begin{equation}
\label{var:LatestDeparture}
    T^{L}_{i} = \left\{
    \begin{aligned}
        & \min\left\{ T^{L}_{i+1} - t_{r_{i}r_{i+1}} - s_{r_{i+1}}, %t\left[r_{i}, r_{i+1}\right],
        l_{r_{i}} + s_{r_{i}} \right\}, \quad \forall i = 1, 2, \dots, n-1 \\
        & |T|\delta, \quad i = n
    \end{aligned}
    \right.
\end{equation}

\begin{equation}
\label{var:SlackForward}
    T^{FS}_{i} =
    \left \{
    \begin{aligned}
        & 0, \quad i=0\\
        & T^{FS}_{i-1} 
        + \max\left\{e_{r_{i}} - \left(T^{E}_{i-1} + s_{r_{i-1}} + t_{r_{i-1}r_i}%t[r_{i-1}, r_{i}]
        \right), 0 \right\}, \quad \text{if } r_{i-1}\in V_c\\
        & \max\left\{e_{r_{i}} - \left(T^{E}_{i-1} + s_{r_{i-1}} + %t\left[r_{i-1}, r_{i}\right]
        t_{r_{i-1}r_i}
        \right), 0 \right\}, \quad \text{otherwise}
    \end{aligned}
    \right.
\end{equation}

\begin{equation}
\label{var:SlackBackward}
    T^{BS}_{i} =
    \left \{
    \begin{aligned}
        & 0, \quad i=n\\
        &  T^{BS}_{i+1}
        + \max\left\{T^{L}_{i+1} 
        - t_{r_{i}r_{i+1}} %t[r_{i}, r_{i+1}] 
        - l_{r_{i}} - s_{r_{i}}, 0 \right\}, 
        \quad \text{if } r_{i+1}\in V_c\\
        &  \max\left\{T^{L}_{i+1} - t_{r_{i}r_{i+1}} %t[r_{i}, r_{i+1}] 
        - l_{r_{i}} - s_{r_{i}} , 
        0 \right\}, \quad \text{otherwise}
    \end{aligned}
    \right.
\end{equation}

\begin{equation}
\label{var:EnerConsump}
    F_i = 
    \left\{
        \begin{aligned}
            & 0, \quad\text{if } i=0 \\
            & F_{i-1} + f_{r_{i-1}r_{i}} %f\left[r_{i-1}, r_{i}\right]
            , \text{ if } r_{i-1}\in V_c \\
            & f_{r_{i-1}r_{i}}
            %f\left[r_{i-1}, r_{i}\right]
            , \text{otherwise.}
        \end{aligned}
    \right.
\end{equation}

$T^{E}_{i}$ is the earliest service start time at node $r_i$ without violating any time window constraints before it. Similarly, $T^{L}_{i}$ is the latest departure time from $r_i$ that will not result in any time window violations after it. 
% Note that, for the calculation of $T^{E}_i$, we assume the arrival time at previous nodes is no later than their corresponding latest service start time for better evaluation of time window violations which is explained in detail later;
$T^{FS}_{i}$ is the forward cumulative slack time, i.e. the difference between earliest arrival time and the earliest service start time, from the last station/depot to $r_i$. $T^{BS}_i$ is the backward cumulative slack time, i.e. the difference between the latest departure time and the latest service time, from $r_i$ to the preceding station/depot.
$F_{i}$ is the amount of time required to recharge the energy consumption from the last station/depot to~$r_i$.

%****************************************
\paragraph{Cargo Capacity Violation}
%****************************************
The cargo capacity violation $\Phi_{cargo}(R)$ of a route $R$ is calculated as
\begin{equation}
\label{func:CargoVio_route}
    \Phi_{cargo}(R) = \max \left\{ \sum_{i = 1}^{n} {q_{r_i}} - Q, 0 \right\}.
\end{equation}

%****************************************
\paragraph{Time Window Violation}
%****************************************

{We calculate the time window violation for route $R$ as the sum of violations at each node along $R$. For the violation at node $r_i$ that follows $r_{i-1}$ with time window violation, we assume the arrival time at $r_{i-1}$ is its latest service start time. This assumption is similar to that of \cite{schneider2014electric} to prevent time window violations from propagating along the route and avoid penalizing a good customer sequence only because it follows a node with time window violation. }The time window violation $\Phi_{tw}(R)$ for a route $R$ is
\begin{equation}
\label{func:TwVio_route}
    \Phi_{tw}(R) = \sum_{i = 1}^{n}\max \left\{ T^{E}_{i} - l_{ r_{i}}, 0 \right\}.
\end{equation}

%****************************************
\paragraph{Battery Capacity Violation}
%****************************************
The battery violations $\Phi_{batt}(R)$ for a route $R$ is
\begin{equation}
\label{func:BattVio_route}
    \Phi_{batt}(R) = \sum_{r \in R \cap V_s} \max\left\{F_r - B, 0\right\}.
\end{equation}

%****************************************
\subsubsection{Net Electricity Cost}
%****************************************

To calculate the electricity cost of a given route, a feasible charging/discharging schedule for each route is computed. Due to the complexity of the problem, the following assumptions are made to limit the number of potential feasible schedules: 
\begin{enumerate}
 \item The maximum number of station visits is limited to two. In practice, it is likely that this assumption is reasonable as the EV fleet needs to primarily service customers, and charging activities are to make sure that the vehicle has enough energy to complete the route and/or to offer the ancillary service of discharging for additional gain. Furthermore, frequent battery charging and discharging have a negative impact on the battery health and as such in practice it is expected that charging and discharging activities are limited to a few station visits.
 \item If a vehicle visits a station to recharge, then the vehicle recharges the battery just enough to be able to complete the route to reach the depot. This assumption is applicable to the two-station visit case and excludes the possibility of over-charging at one station to discharge at another station later.
 \item {A vehicle is allowed to perform either charging or discharging at stations, but is not allowed to do both during a single station visit. This assumption reduces the solution space defined by the MIP formulation.}% but the solution quality is not significantly compromised according to the computational experiment shown in Section \ref{compresults}.
\end{enumerate}

Given a depot or a station node which is the $i^{th}$ node on route $R$ denoted by $r_i \in R$, the time window during which a vehicle can charge or discharge is $\left[T^{E}_{r_i}, T^{L}_{r_i}\right]$. Given that the planning horizon is discretized, then the set of time periods during which the vehicle can charge/discharge is given by $T_i = \left\{\lceil\frac{T^{E}_{r_i}}{\delta}\rceil, \lceil\frac{T^{E}_{r_i}}{\delta}\rceil + 1, \dots, \lfloor\frac{T^{L}_{r_i}}{\delta}\rfloor\right\}$, namely the \textit{connected periods} for $r_i$.

Given the $i^{th}$ and the $j^{th}$ node ($i<j$) of route $R$ where a vehicle can charge/discharge (station or depot node) and $T_i$ and $T_j$ are the connected periods corresponding to $i$ and $j$ respectively, if a vehicle charges/discharges at time period $t \in T_i$, then this might make it impossible for the vehicle to arrive at node $j$ before the start of time period $t'$. Thus $M^t_i$ is defined as the \textit{mutually exclusive set} for $t\in T_i$ which is the set of all the time periods $t'\in T_j$ where the vehicle cannot reach node $j$ if it is charging/discharging during time period $t$. The mutually exclusive set $M^t_i$ is given by 
\begin{equation}
    \label{var:MutuallyExclusiveSet}
    M_i^t = 
    \left\{
    \begin{aligned}
        & \emptyset,
        \text{ if } \delta t \leq T^{E}_{r_i} + T^{FS}_{r_j}\\
        & \left\{t': \delta (t'-1) 
        < T^{E}_{r_j} + \delta t -\left( T^{E}_{r_i} + T^{FS}_{r_j}\right)
        \right\},
        \text{ otherwise.}
    \end{aligned}
    \right.
\end{equation}

Given that at most two stations are allowed in each route, three different cases are analyzed to find a feasible charging/discharging schedule corresponding to the three special cases of zero stations, one station, and two stations routes, respectively.
\paragraph{Zero stations} If there are no stations along the vehicle's route and given that the EV is fully charged at the starting depot, then the only potential activities are to discharge the EV at the starting and/or ending depot. The maximum amount of time periods during which the EV can discharge while ensuring that there is enough energy available to cover the full route is given by $\Omega = \lfloor \frac{B - F_n}{\delta} \rfloor$. Given the connected periods  $T_1$ and $T_n$ at the two depot nodes (the first and the last node along $R$) and the set of mutually exclusive time periods $M_1^t$ for each time period $t\in T_1$, the following mixed integer program maximizes the revenue from energy discharge
\begin{align}
\label{obj:NoStat}
    \Upsilon^* = \max_d & \sum_{t \in T_1}P_{dis}^{t}d^t_{1} + \sum_{t\in T_n}P_{dis}^{t}d_{n}^{t}\\
\label{con:NoStat_BattCap}
    \textrm{s.t.} \quad & \sum_{t\in T_1}d^t_{1} + \sum_{t\in T_n}d^t_{n} \leq \Omega,\\
\label{con:NoStat_MutuallyExc}
    & d^{t}_{1} + d^{k}_{n} \leq 1, \quad \forall t \in T_1, \ \forall k \in M^t_1, \\
\label{con:NoStat_domaim}
    & d^t_j \in \{0, 1\}, \quad \forall j\in \{1, n\},\ \forall t\in T_j.
\end{align}
where $d^m_j$ is a binary decision variable indicating if the EV discharges at node $j$ during time period~$m$. We note that by enumerating the end period where the EV discharges at $r_1$, the MIP simplifies to $|T_1|$ $0$-$1$ knapsack problems with unit item weights that can be solved very efficiently. Given the optimal solution $d^*$, the net cost of electricity for route $R$ is given by 
\begin{equation}
    \label{func:ElecCost_NoStat}
    f_{elec}(R) = - \Upsilon^* + (F_n + \sum_{t \in T_1}{\delta d_1^{t}}^* + \sum_{k\in T_n} {\delta d_n^{k}}^*) \times P_{night}.
\end{equation}

\paragraph{One station} Given that an EV is visiting a station ($k^{th}$ node along route $R$), then one possibility is that the EV has to recharge the battery at the station node to make sure there is enough energy to reach the end-of-route depot ($n^{th}$ node along route $R$).
This case occurs if $\Delta = F_k + F_n - B  > 0$. Otherwise, the EV can either discharge energy at the station or recharge in order to discharge later at the depot. Given the sets of connected time periods $T_1$, $T_k$, and $T_n$, for the $1^{st}$, $k^{th}$, and last node along route $R$, respectively, and the mutually exclusive sets $M^t_1$ for each time period $t\in T_1$, and $M^t_k$ for each time period $t\in T_k$, the following mixed integer program maximizes the revenue from energy charging and discharging
\begin{align}
\label{obj:OneStat}
    \Upsilon^* =  \max_{d, r} \ &  \sum_{t \in T_1}P_{dis}^{t}d^t_{1} + \sum_{t\in T_k}P_{dis}^{t}d^t_{k} + \sum_{t\in T_n}P_{dis}^{t}d^t_{n} -  \sum_{t\in T_k}P_{re}^{t}r^t_{k},\\
\label{con:OneStat_BattCapDown_Dis}
    \mbox{s.t. }& \sum_{t \in T_1}d^t_{1} + \sum_{t\in T_k}d^t_{k} + \sum_{t\in T_n}d^t_{n} - \sum_{t\in T_k}r^t_{k}
    \begin{cases}
     \leq \lfloor \frac{-\Delta}{\delta} \rfloor &\mbox{ if $\Delta < 0$,}\\
     = \lceil \frac{\Delta}{\delta} \rceil &\mbox{ otherwise,}
    \end{cases}
 \\
\label{con:OneStat_BattCapUp}
    &\delta \sum_{t\in T_k}r^t_{k}  \leq F_k + \delta \sum_{t \in T_1}d^t_{1},\\
\label{con:OneStat_ArrDepot}
    &\delta\left(\sum_{t \in T_1}d^t_{1} + \sum_{t\in T_k}d^t_{k}\right)  + F_k + F_n  \leq B + \delta \sum_{t\in T_k}r^t_{k},\\
\label{con:OneStat_Dis}
    &\sum_{t\in T_k}d^t_{k} \leq |T_k|y,\\
\label{con:OneStat_Re}
     &\sum_{t\in T_k}r^t_{k} \leq |T_k| (1-y),\\
\label{con:OneStat_MutuallyExc1}
    &d^{t_1}_1 + r^{t_2}_k + d^{t_2}_k \leq 1, \quad \forall t_1 \in T_1, \ \forall t_2 \in M^{t_1}_1,\\
\label{con:OneStat_MutuallyExc2}
    &r^{t_1}_k + d^{t_1}_k + d^{t_2}_n \leq 1, \quad \forall t_1 \in T_k, \ \forall t_2 \in M^{t_1}_k, \\
\label{con:OneStat_domain_rd}
    & d^t_j, r^t_j \in \{0, 1\}, \quad \forall j\in \{1, k\},\ \forall t\in T_j,\\
\label{con:OneStat_domain_y}
    & y \in \{0, 1\}.
\end{align}
{where $d^m_j$ and $r^m_j$ are binary decision variables indicating if the EV discharges and charges at node $j$ during time period~$m$, respectively. Constraint \eqref{con:OneStat_BattCapDown_Dis} specifies the amount of charging and discharging given different values of $\Delta$. Constraint \eqref{con:OneStat_BattCapUp} ensures the EV battery capacity is not violated after charging/discharging at the en-route station. Constraint \eqref{con:OneStat_ArrDepot} guarantee the EV has enough energy to complete the trip. Constraints \eqref{con:OneStat_Dis}-\eqref{con:OneStat_Re} force that the EV can perform either charging or discharging at the station. Constraints \eqref{con:OneStat_MutuallyExc1}-\eqref{con:OneStat_MutuallyExc2} describe the mutual exclusiveness among periods at $r_1$, $r_k$ and $r_n$. Constraint \eqref{con:OneStat_domain_rd}-\eqref{con:OneStat_domain_y} specify the domains of the decision variables.} Given the optimal solution $d^*$ and $r^*$, the net cost of electricity for route $R$ is given by
\begin{equation}
\label{func:ElecCost_OneStat}
    f_{elec}(R) = -\Upsilon^* + (F_k + F_n + \sum_{i \in T_1}{\delta d_1^{i}}^* + \sum_{j \in T_k}{\delta d_k^{j}}^* + \sum_{t\in T_n} {\delta d_n^{t}}^* - \sum_{j\in T_k}{\delta r^j_{k}}^* ) \times P_{night}.
\end{equation}

\paragraph{Two stations} Given an EV route that includes visits to two stations $k_1$ and $k_2$ with $k_1<k_2$. If $\Delta=F_{k_1} + F_{k_2} +F_{n} - B>0$ then the EV needs to recharge en-route in order to complete the trip to the final depot. In order to formulate the problem to minimize the cost of energy recharging, we first define the sets of connected time periods $T_{k_1}$ and $T_{k_2}$ corresponding to stations $k_1$ and $k_2$ respectively. Furthermore, let $M^t_{k_1}$ be the mutually exclusive sets for each time period $t\in T_{k_1}$. Given the binary variables $r^t_{k_j}$,
%\begin{align*}
%r^t_{k_j} = \begin{cases}
%           1 \mbox{ if EV recharges at %station $k_1$ during time period t},\\
%           0 \mbox{ otherwise,}
%         \end{cases}
%\end{align*}
the problem that minimizes the cost of energy recharging is given by
\begin{align}
\label{obj:TwoStat_re}
    \Upsilon^* =  \min_{r} \ & \sum_{t\in T_{k_1}}P_{re}^{t}r^t_{k_1} + \sum_{t \in T_{k_2}}P_{re}^{t}r^t_{k_2},\\
\label{con:TwoStat_re_BattCapAll}
    \mbox{s.t. }&\quad \sum_{t\in T_{k_1}}r^t_{k_1} + \sum_{t \in T_{k_2}}r^t_{k_2} = \Omega,\\
\label{con:TwoStat_re_k1}
     &F_{k_1} + F_{k_2} - B \leq \delta \sum_{t\in T_{k_1}}r^t_{k_1} \leq F_{k_1},\\
\label{con:TwoStat_re_k2}
    &F_{k_2} + F_{n} - B \leq \delta \sum_{t\in T_{k_2}}r^t_{k_2} \leq F_{k_1} + F_{k_2} - \delta\sum_{t \in T_{k_1}}r^t_{k_1},\\
\label{con:TwoStat_re_MutuallyExc}
    &r^{t_1}_{k_1} + r^{t_2}_{k_2} \leq 1, \quad \forall t_1 \in T_{k_1}, \ \forall t_2 \in M^{t_1}_{k_1}, \\
\label{con:TwoStat_domain}
    & r^t_j \in \{0, 1\}, \quad \forall j\in \{k_1, k_2\}, \ \forall t\in T_j.
\end{align}
The objective function \eqref{obj:TwoStat_re} minimizes the en-route electricity recharging cost. Constraints \eqref{con:TwoStat_re_BattCapAll} make sure the EV has enough energy to complete the trip. Constraints \eqref{con:TwoStat_re_k1} ensure that the amount of recharged energy at the first station does not exceed the available battery capacity and will allow the EV to reach the next station. Similarly, constraints \eqref{con:TwoStat_re_k2} enforce the capacity limits at the second station. Finally, constraints \eqref{con:TwoStat_re_MutuallyExc} enforce the mutual exclusive conditions for the time periods. Given $r^*$, the optimal solution of problem \eqref{obj:TwoStat_re}--\eqref{con:TwoStat_re_MutuallyExc}, the route's cost is then given by
\begin{equation}
\label{func:TwoStat_ElecCost_re}
    f_{elec}(R) = \Upsilon^* + \left(F_{k_1}+F_{k_2}+F_{n} - \sum_{t\in T_{k_1}}\delta r^t_{k_1} - \sum_{t \in T_{k_2}}\delta r^t_{k_2}\right) \times P_{night}.
\end{equation}

If $\Delta<0$, then the EV can discharge energy at the depots (nodes 1 and $n$) or the stations (nodes $k_1$ and $k_2$). To formulate the problem that maximizes the value of the discharged energy, we introduce binary decision variable $d^{t}_{k}$ that takes value 1 if the EV discharges at node $k$ during period $t$ and takes value $0$ otherwise.
% \begin{align*}
% d^t_{k} = \begin{cases}
%           1 \mbox{ if EV discharges energy at node $k$ during time period $t$},\\
%           0 \mbox{ otherwise.}
%           \end{cases}
% \end{align*}
Given the connected time periods $T_1$, $T_{k_1}$, $T_{k_2}$, and $T_{n}$ corresponding to the starting depot, stations $k_1$ and $k_2$, and the ending depot, respectively, and $M^t_1$, $M^t_{k_1}$, and $M^t_{k_2}$ the mutually exclusive sets for each time period $t\in T_1$, $t\in T_{k_1}$, and $t\in T_{k_2}$, the problem is formulated as
\begin{align}
\label{obj:TwoStat_dis}
    \Upsilon^* =  \max_{d} \ & \sum_{t\in T_{1}}P_{dis}^{t}d^t_{1} + \sum_{t\in T_{k_1}}P_{dis}^{t}d^t_{k_1} + \sum_{t \in T_{k_2}}P_{dis}^{t}d^t_{k_2} + \sum_{t \in T_{n}}P_{dis}^{t}d^t_{n}\\
\label{con:TwoStat_dis_BattCapAll}
    \mbox{s.t. }& \sum_{t\in T_{1}}d^t_{1} + \sum_{t\in T_{k_1}}d^t_{k_1} + \sum_{t \in T_{k_2}}d^t_{k_2} + \sum_{t \in T_{n}}d^t_{n} \leq \Omega\\
\label{con:TwoStat_dis_MutuallyExc}
    & d^{t_1}_{1} + d^{t_2}_{k_1} \leq 1, \quad \forall {t_1} \in T_{1}, \ \forall {t_2} \in M_{1}^{t_1},\\
    & d^{t_1}_{k_1} + d^{t_2}_{k_2} \leq 1, \quad \forall {t_1} \in T_{k_1}, \ \forall {t_2} \in M_{k_1}^{t_1},\\
    & d^{t_1}_{k_2} + d^{t_2}_{n} \leq 1, \quad \forall {t_1} \in T_{k_2}, \ \forall {t_2} \in M_{k_2}^{t_1}.
\end{align}
 Given the optimal solution $d^*$, the net cost of electricity for route $R$ is given by
\begin{equation}
\label{func:TwoStat_ElecCost_dis}
    f_{elec}(R) = - \Upsilon^{*} + \left(F_{k_1}+F_{k_2}+F_{n} + \sum_{t\in T_{1}} \delta  d^t_{1} + \sum_{t\in T_{k_1}} \delta d^t_{k_1} + \sum_{t \in T_{k_2}} \delta d^t_{k_2} + \sum_{t \in T_{n}} \delta  d^t_{n}\right) \times P_{night}.
\end{equation}

%======================================================================
\subsection{Variable Neighborhood Search}
%======================================================================

\begin{figure}[ht]
    \centering
    \includegraphics[width = 5in]{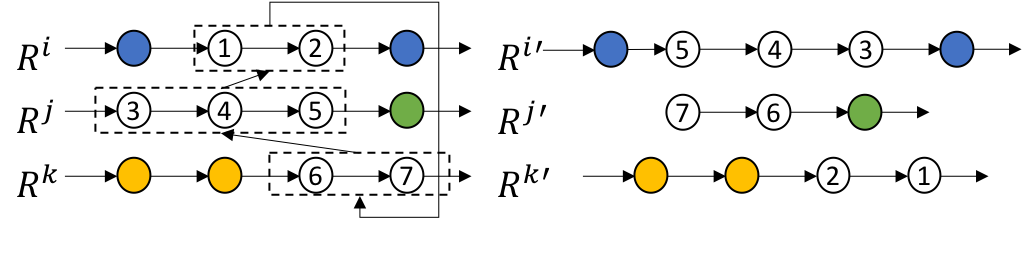}
    \caption{An Example of the Cyclic Exchange Operator}
    \label{fig:CyclicExchange}
\end{figure}

The variable neighborhood search heuristic was proposed by \cite{mladenovic1997variable}. Given the current solution $S$, the neighborhood structure is defined by a cyclic-exchange operator which selects $N_{r}$ routes from $S$ to form an exchange cycle. The cyclic-exchange selects in each route $R^i$, a random number $\upsilon_i$ of consecutive nodes which form an exchange block. These blocks are then reversed and exchanged between the routes. Each exchange of blocks forms a new neighboring solution $S^{'}$ which may be feasible or infeasible.
Figure \ref{fig:CyclicExchange} shows an example of the cyclic exchange operator where $N_r = 3$. The selected routes are $R^i$, $R^j$ and $R^k$ with $\upsilon_i = 2$, $\upsilon_j = 3$, $\upsilon_k = 2$, respectively. 
The three blocks on the left form an exchange cycle. 
The blocks are reversed and transferred forming the three new routes ${R^{i}}'$, ${R^{j}}'$, and ${R^{k}}'$ on the right.

Given a neighboring solution $S^{'}$, instead of the commonly used approach of applying local descent to improve the solution, a tabu search is applied similar to \cite{schneider2014electric} to find a local optima $S^{''}$ which is accepted if it is better than the current solution $S$. To further diversify the search, $S^{''}$ is also accepted with a probability of $e^{\left[\frac{f_{gen}(S) - f_{gen}(S^{''})}{Temp}\right]}$ if it is worse than the current solution $S$ \citep{schneider2014electric, hemmelmayr2009variable, stenger2013adaptive}. The temperature $Temp$ of this annealing phase is initialized to $Temp_0$ such that a solution with cost $f_{gen}(S^{''})$ that is $\kappa$ times worse than the current best solution will be accepted with a probability of 50\%. The temperature is then linearly decreased by a factor after each VNS iteration so that in the last 20\% of the iterations, the temperature is below $0.0001$. These parameters are the same as the ones used by \cite{schneider2014electric}.

%======================================================================
\subsection{Tabu Component}
%======================================================================

\begin{figure}[!ht]
    \centering
    \includegraphics[width = 4.5in]{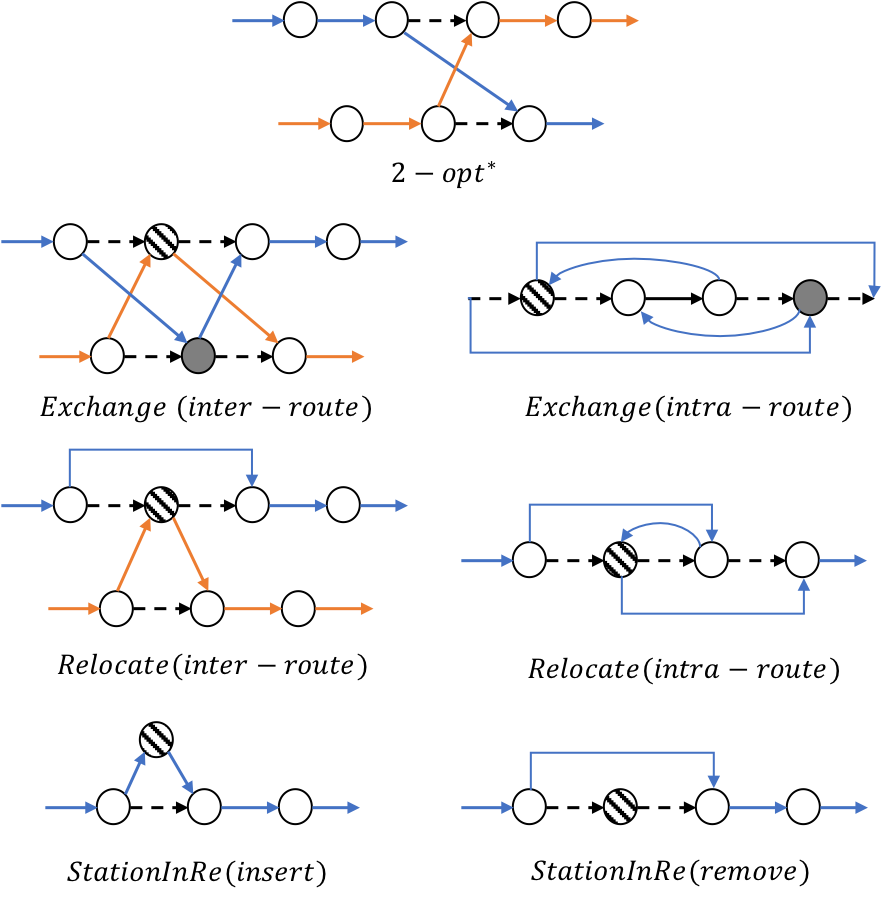}
    \caption{Tabu Search Operators}
    \label{fig:TabuOperators}
\end{figure}

Tabu search (TS) is applied to every solution $S^{'}$ generated by VNS. The neighborhood solutions are generated using the three widely-used operators for tabu search: 2-opt$^*$, exchange, and relocate, as well as the stationInRe operator that is discussed in \cite{schneider2014electric}. The operators are visualized in Figure \ref{fig:TabuOperators} where the nodes that are on the same row are travelled by an EV before the operator is applied, the dashed arrows are edges to be removed, and the stripped and shadowed nodes are the nodes that are selected by the algorithm. The operators are applied as follows:
\begin{itemize}
    \item \textbf{2-opt*:} Select two routes and remove one edge from each of them. Connect the first part of the first route with the second part of the second route and vice versa.
    \item \textbf{Exchange:} Exchange the positions of two nodes. The two nodes could either be in the same route or in different routes.
    \item \textbf{Relocate:} Select one route and remove one node from this route. Reinsert the selected node at another position. The new position could either be in the same route or in another route.
    \item \textbf{StationInRe:} Perform insertion or removal of a station node.
\end{itemize}

At each TS iteration, the four operators are applied on $S^{'}$ and the candidate solutions are filtered against the tabu list. The solution that results in the maximal decrease in the generalized cost $f_{gen}(S)$ is selected. The reinsertion of the removed edges which led to the selected solution, is prohibited for a fixed number of tabu iterations called tabu tenure. The tabu tenure for each deleted edge is randomly selected from an interval $\left[v_{min}, v_{max} \right]$. The procedure is repeated until no improvement can be made or for at most $\eta_{tabu}$ iterations. The best solution that is obtained is denoted as $S^{''}$. 

%======================================================================
\section{Computational Experiments} \label{compresults}
%======================================================================

This section presents extensive computational results to evaluate the proposed Lagrangian relaxation approach and the VNS/TS hybrid heuristic.

%======================================================================
\subsection{Test Instances} \label{sec:TestInstance}
%======================================================================

The test instances that are used in the evaluation are based on the ones developed by \cite{schneider2014electric} which were constructed using the instances proposed by \cite{solomon1987algorithms}. % Particularly, we consider instances with $5$, $10$, and $15$ customers with up to $6$ stations which we refer to as the Schneider instances.
The instances are classified into $3$ categories based on the geographical distribution of the customer nodes (see first column of Table \ref{BoundEvaluation}). The instances that start with ``R'' are random instances where customers are uniformly distributed, and
those that start with ``C'' are clustered instances in which customers are clustered into small groups. The customer distribution of the ``RC'' instances is a mixture of random and clustered distributions. 

{The EV fleet is homogeneous, each with a cargo capacity $Q = 200$, a battery capacity $B = 270$, and a range of $150$ km. The travel speed $v$ is set to a constant value of $30$ km/h, and the batteries are recharged at a constant speed $\alpha = 1.8$. The discharging speed is assumed the same as the charging speed. The choice of these values are discussed in more details in the case study.}

{The locations and time windows in the Schneider instances are normalized values which makes it difficult to relate the charging and discharging time to real-world values so as to estimate the potential costs and gains. We adjust the instances in a way similar to \citet{schiffer2018strategic}. In particular, we set the maximum distance from the depot to customers/stations as $100$ km and convert all other distances to km proportionally for each instances. In addition, since fully charging an EV requires $B = 270$ minutes ($4.5$ hours), we set the planning horizon as $5am - 12am$ ($\left[0, 1140\right]$ in minutes) so that an EV has enough time to recharge at night. We scale all time windows proportionally to minutes to fit the planning horizon. We find through computational experiments that the time windows in the scaled instances are relatively tight which prevent vehicles from detouring to perform charging and discharging activities. We thus further relax the time windows to three periods, morning ($5am - 12pm$), afternoon ($12pm - 6pm$), and evening ($6pm - 12am$), based on their service start time. }

For charging and discharging, the length of each period is set to one hour, i.e. $\delta = 60$ minutes. The energy prices which are shown in Table \ref{tab:hydro_rate} are based on the real time-of-use hydro rate in Ontario, Canada in effect between May $1$, $2019$ and October $31$, $2019$. The reward rates are chosen such that an EV can make profits by discharging at peak hours (11:00 AM - 5:00 PM) and recharging the battery later. Note that the reward rates are also economically beneficial to the grid because the discharging reward that the grid pays to the EV owner is lower than the corresponding market price. 

\begin{table}
\centering
\caption{Time-of-Use Electricity Prices}
\begin{tabular}{@{}cccc@{}}
\toprule
From     & To       & Charging Cost ($\cent$/kWh) & Discharging Reward ($\cent$/kWh) \\ \midrule
12:00 AM & 7:00 AM  & 6.5           & 6.5                \\
7:00 AM  & 11:00 AM & 9.4           & 8.0                  \\
11:00 AM & 5:00 PM  & 13.4          & 10.0                 \\
5:00 PM  & 7:00 PM  & 9.4           & 8.0                  \\
7:00 PM  & 12:00 AM & 6.5           & 6.5                \\ \bottomrule
\end{tabular}
\label{tab:hydro_rate}
\end{table}

%======================================================================
\subsection{Experimental Setup} \label{sec:TestSettings}
%======================================================================

All the tests are performed on a MacBook Pro running OS X 10.13.6, using a single $2.30$ GHZ CPU, $16$ GB of RAM and CPLEX $20.1.0.0$ is used as an optimization solver. The time limit is set to $7200$ seconds and the memory limit is set to $10$ GB. The algorithms are implemented as single thread codes in Python.

For the VNS/TS hybrid heuristic, the penalty parameters $\beta_{tw}$, $\beta_{batt}$, $\beta_{cargo}$ are set to $10$ and the number of tabu iterations per round is set as $\eta_{tabu} = 30$. The number of VNS iterations $\eta_{vns}$ and the early stopping criterion $\eta_{early}$ are varied with respect to the number of customers included. For the instances with 5, 10, and 15 customers, $\eta_{vns}$ is set to $10$, $20$, and $30$ respectively. $\eta_{early}$ is set to 10 for instances with more than 10 customers and 5 otherwise.
For the cyclic operator, $N_r$ is set to $2$ when the fleet consists of $3$ EVs or less and to $3$ otherwise.
The length of each exchange block is randomly selected from $\{1, 2, 3\}$.
The upper and lower bounds of the tabu tenure are $v_{min} = 5$ and $v_{max} = 15$. Parameter $\kappa$ for the annealing mechanism is set to $0.5$.

%======================================================================
\subsection{Computational Results} \label{sec:performance}
%======================================================================

The results for CPLEX, the VNS/TS heuristic, the Lagrangian relaxation and the LP relaxation are presented in Table \ref{BoundEvaluation}. 
The upper and lower bounds achieved by the four algorithms are presented in the ``UB'' and ``LB'' columns, respectively. The ``Time'' column indicates the computational time in seconds. 
The ``Gap'' columns indicate the gap between the corresponding upper bound and the lower bound obtained by the Lagrangian relaxation. {Column ``Best Iter'' indicates the VNS iteration during which the best solution was found.}
\begin{table}[!h]
\caption{Performance of CPLEX, Lagrangian relaxation, and VNS/TS heuristic on small instances}
\centering
\begin{threeparttable}
\scriptsize
{
\begin{tabular}{@{}cclcrrlcrcclcrlcr@{}}
\toprule
\multicolumn{2}{c}{Instance} &  & \multicolumn{3}{c}{CPLEX}  &  & \multicolumn{4}{c}{Heuristic} &  & \multicolumn{2}{c}{Lagrangian} &  & \multicolumn{2}{c}{LP Relaxation}
\\ \cmidrule(r){1-2} \cmidrule(lr){4-6} \cmidrule(lr){8-11} \cmidrule(l){13-14} \cmidrule(lr){16-17} 
Name             & K         &  & UB      & Time    & Gap     &  & UB      & Time  & best iter & Gap    &  & LB             & Time                  &  & LB       & Time\\ \midrule
C101-5           & 2         &  & 481.35  & 2.05    & 0.00\%  &  & 481.35  & 5.58  & 9         & 0.00\% &  & 481.35         & 5.13               &  & -2041.15 & 0.02         \\
C103-5           & 2         &  & 316.83  & 1.50    & 0.00\%  &  & 316.83  & 9.16  & 1         & 0.00\% &  & 316.83         & 3.39               &  & -1522.90 & 0.03         \\
C206-5           & 2         &  & 553.42  & 87.35   & 0.00\%  &  & 553.42  & 10.26 & 1         & 0.00\% &  & 553.42         & 33.23              &  & -2533.86 & 0.03         \\
C208-5           & 2         &  & 471.10  & 22.58   & 0.00\%  &  & 471.10  & 9.40  & 1         & 0.00\% &  & 471.10         & 17.62              &  & -2040.28 & 0.02         \\
R104-5           & 2         &  & 334.86  & 21.34   & 0.00\%  &  & 349.26  & 7.33  & 1         & 4.30\% &  & 334.86         & 8.03               &  & -2041.51 & 0.02         \\
R105-5           & 2         &  & 405.05  & 16.06   & 0.00\%  &  & 409.37  & 12.51 & 4         & 1.07\% &  & 405.05         & 5.07               &  & -2042.08 & 0.02         \\
R202-5           & 1         &  & 530.25  & 1.13    & 0.00\%  &  & 530.25  & 1.85  & 2         & 0.00\% &  & 530.25         & 7.56               &  & -1085.40 & 0.03         \\
R203-5           & 2         &  & 577.86  & 251.03  & 0.00\%  &  & 577.86  & 15.90 & 1         & 0.00\% &  & 577.86         & 101.93             &  & -2519,02 & 0.03         \\
RC105-5          & 2         &  & 482.45  & 261.85  & 0.00\%  &  & 482.45  & 11.42 & 1         & 0.00\% &  & 482.45         & 103.90             &  & -2520.12 & 0.03         \\
RC108-5          & 2         &  & 616.17  & 57.19   & 0.00\%  &  & 616.17  & 10.82 & 1         & 0.00\% &  & 616.17         & 34.30              &  & -2507.49 & 0.03         \\
RC204-5          & 1         &  & 563.02  & 32.99   & 0.00\%  &  & 579.97  & 32.52 & 3         & 3.01\% &  & 563.02         & 129.35             &  & -1409.40 & 0.03         \\
RC208-5          & 1         &  & 465.64  & 0.56    & 0.00\%  &  & 465.64  & 1.41  & 6         & 0.00\% &  & 465.64         & 7.75              &  & -1085.40 & 0.03         \\
                 &           &  &                 &         & &  &        &        &           &        &  &               &                      &  &  \\
C101-10          & 3         &  & $1046.11^{(*)}$ & 4565.51 & 32.49\% &  & 789.58 & 62.68 & 10         & 0.00\% &  & 789.58        & 2194.74        &  & -4802.90 & 0.09\\
C104-10          & 2         &  & $867.31^{(*)}$  & 952.55  & -       &  & 702.52 & 54.62 & 2          & -      &  & -             & $>$7200.00     &  & -2818.80 & 0.07\\
C202-10          & 2         &  & $700.99^{(*)}$  & 3766.34 & 15.39\% &  & 607.49 & 91.60 & 9          & 0.00\% &  & 607.49        & $>$7200.00     &  & -3466.80 & 0.07\\
C205-10          & 3         &  & 616.36          & 4071.78 & 1.88\%  &  & 616.36 & 95.15 & 3          & 1.88\% &  & 604.98        & 842.72         &  & -3256.20 & 0.05\\
R102-10          & 4         &  & $1090.31^{(*)}$ & 1976.55 & 96.01\% &  & 561.96 & 89.80 & 8          & 1.02\% &  & 556.26        & $>$7200.00     &  & -5022.62 & 0.08\\
R103-10          & 2         &  & $458.02^{(*)}$  & 1463.56 & -       &  & 443.75 & 55.05 & 1          & -      &  & -             & $>$7200.00     &  & -2170.80 & 0.03\\
R201-10          & 2         &  & $713.81^{(*)}$  & 1752.25 & 20.81\% &  & 590.83 & 46.31 & 1          & 0.00\% &  & 590.83        & $>$7200.00     &  & -2818.80 & 0.06\\
R203-10          & 2         &  & -               & -       & -       &  & 742.42 & 91.92 & 17         & -      &  & -             & $>$7200.00     &  & -3466.80 & 0.07\\
RC102-10         & 4         &  & $1211.21^{(*)}$ & 1558.34 & 26.58\% &  & 999.48 & 80.09 & 2          & 4.45\% &  & 956.87        & 1931.84        &  & -5037.47 & 0.10\\
RC108-10         & 4         &  & $1221.72^{(*)}$ & 1350.84 & 57.53\% &  & 775.57 & 135.29& 12         & 0.00\% &  & 775.57        & $>$7200.00     &  & -5018.51 & 0.09\\
RC201-10         & 2         &  & $850.35^{(*)}$  & 2041.85 & 12.02\% &  & 761.15 & 81.25 & 1          & 0.27\% &  & 759.08        & 1850.80        &  & -2818.80 & 0.05\\
RC205-10         & 3         &  & $1243.35^{(*)}$ & 1448.82 & 46.07\% &  & 851.21 & 73.49 & 2          & 0.00\% &  & 851.21        & $>$7200.00     &  & -4020.07 & 0.07\\
                 &           &  &        &         &         &  &               &                &  &        &        &           &        \\
C103-15          & 3         &  & *              & 2102.03   & -       &  & 770.22   & 528.67 & 14         & -   &  & -             & $>$7200.00              &  & -5200.20 & 0.15\\
C106-15          & 3         &  & $670.87^{(-)}$ & $>$7200.00& -       &  & 657.32   & 320.45 & 10        & -    &  & -             & $>$7200.00              &  & -3256.20 & 0.08\\
C202-15          & 2         &  & -              & 3303.37   & -       &  & 1019.80  & 271.92 & 3        & -     &  & -             & $>$7200.00              &  & -3466.80 & 0.10 \\
C208-15          & 2         &  & *              & 2875.61   & -       &  & 757.12   & 139.56 & 2         & -    &  & -             & $>$7200.00              &  & -2818.80 & 0.07  \\
R102-15          & 5         &  & *              & 1975.12   & -       &  & 692.38   & 407.32 & 5         & -    &  & -             & $>$7200.00              &  & -12295.33 & 0.35  \\
R105-15          & 4         &  & *              & 1722.53   & -       &  & 790.48   & 392.96 & 4        & -     &  & -             & $>$7200.00              &  & -7638.94 & 0.17 \\
R202-15          & 3         &  & -              & $>$7200.00& -       &  & 1114.25  & 299.90 & 3         & -    &  & -             & $>$7200.00              &  & -5981.33 & 0.07  \\
RC103-15         & 4         &  & *              & 2454.94   & -       &  & 667.52   & 495.01 & 14        & -    &  & -             & $>$7200.00              &  & -6505.45 & 0.17    \\
RC108-15         & 3         &  & *              & 1660.90   & -       &  & 873.22   & 206.70 & 1        & -     &  & -             & $>$7200.00              &  & -5200.20 & 0.12 \\
RC202-15         & 3         &  & *              & 2909.87   & -       &  & 996.73   & 267.89 & 2         & -    &  & -             & $>$7200.00              &  & -5200.20 & 0.21  \\
RC204-15         & 3         &  & -              & $>$7200   & -       &  & 971.64   & 270.05 & 3         & -    &  & -             & $>$7200.00             &  & -4762.80 & 0.15  \\ \bottomrule
\end{tabular}
}
\begin{tablenotes}
    \footnotesize
    \item Instances that violate the memory and time limits are labeled with $*$ and $-$ respectively
\end{tablenotes}
\end{threeparttable}
\label{BoundEvaluation}
\end{table}

\begin{table}[!h]
\caption{Performance of VNS/TS heuristic on medium instances}
\centering
\scriptsize
{
\begin{tabular}{@{}ccccccccccccc@{}}
\toprule
\multicolumn{2}{c}{Instance} &  & \multicolumn{3}{c}{Heuristic} &  & \multicolumn{2}{c}{Instance} &  & \multicolumn{3}{c}{Heuristic} \\ \cmidrule(r){1-2} \cmidrule(lr){4-6} \cmidrule(lr){8-9} \cmidrule(l){11-13} 
Name             & K         &  & UB      & Time    & Best iter &  & Name             & K         &  & UB      & Time     & Best iter \\ \midrule
C103-20          & 4         &  & 1049.71 & 1147.21 & 20        &  & C103-30          & 4         &  & 1059.30  & 3340.14  & 2         \\
C106-20          & 3         &  & 833.62  & 495.04  & 2         &  & C106-30          & 4         &  & 1724.82  & 3951.79  & 7        \\
C202-20          & 4         &  & 1202.31 & 642.04  & 1         &  & C202-30          & 4         &  & 1477.16  & 3615.39  & 5        \\
C208-20          & 3         &  & 848.83  & 385.78  & 1         &  & C208-30          & 4         &  & 1238.15  & 2798.32  & 3        \\
R102-20          & 3         &  & 1385.19 & 713.47  & 2         &  & R102-30          & 4         &  & 1423.15  & 3036.66  & 3         \\
R105-20          & 4         &  & 1257.00 & 789.36  & 2         &  & R105-30          & 4         &  & 1895.08  & 2883.65  & 4        \\
R202-20          & 3         &  & 1096.58 & 642.69  & 3         &  & R202-30          & 4         &  & 1860.16  & 2447.00  & 4         \\
RC103-20         & 3         &  & 1152.61 & 790.07  & 17        &  & RC103-30         & 4         &  & 1674.06  & 3182.44  & 17        \\
RC108-20         & 4         &  & 1358.33 & 944.22  & 6         &  & RC108-30         & 4         &  & 1579.86  & 3998.71  & 17         \\
RC202-20         & 3         &  & 1545.39 & 1101.57 & 11        &  & RC202-30         & 5         &  & 1850.08  & 3608.25  & 9        \\
RC204-20         & 3         &  & 996.12  & 745.32  & 5         &  & RC204-30         & 4         &  & 1391.72  & 2297.52  & 1        \\
                 &           &  &         &         &           &  &                  &           &  &         &          &           \\
C103-25          & 4         &  & 857.06  & 2191.81 & 15        &  & C103-35          & 4         &  & 1324.19  & 5447.90  & 11         \\
C106-25          & 3         &  & 1329.92 & 1262.92 & 11        &  & C106-35          & 4         &  & 1152.76  & 4954.01  & 12        \\
C202-25          & 3         &  & 1154.42 & 1155.89 & 1         &  & C202-35          & 4         &  & 1759.25  & 2801.60  & 2       \\
C208-25          & 3         &  & 1082.15 & 796.42  & 2         &  & C208-35          & 4         &  & 1520.57  & 3924.57  & 11       \\
R102-25          & 4         &  & 1164.46 & 1495.67 & 1         &  & R102-35          & 4         &  & 1928.67  & 3851.51  & 8         \\
R105-25          & 4         &  & 1811.40 & 1088.10 & 4         &  & R105-35          & 4         &  & 1639.93  & 2854.46  & 6        \\
R202-25          & 3         &  & 1209.24 & 1632.27 & 17        &  & R202-35          & 4         &  & 1412.07  & 4365.94  & 14        \\
RC103-25         & 4         &  & 1372.35 & 2573.81 & 16        &  & RC103-35         & 5         &  & 1857.54  & 3135.70  & 1        \\
RC108-25         & 4         &  & 1341.98 & 2416.85 & 10        &  & RC108-35         & 5         &  & 1786.30  & 3906.11  & 2         \\
RC202-25         & 4         &  & 1472.85 & 2173.26 & 7         &  & RC202-35         & 4         &  & 1920.45  & 3047.78  & 3        \\
RC204-25         & 4         &  & 1644.16 & 1896.30 & 3         &  & RC204-35         & 5         &  & 1268.05  & 6068.46  & 17        \\ 
\bottomrule
\end{tabular}
}
\label{BoundEvaluation_medium_instance}
\end{table}

For the small instances with $5$ customers, the VNS/TS heuristic outperforms CPLEX in terms of solution time for the majority of the tested instances. In terms of bound quality, the Lagrangian relaxation obtains the optimal solution for every instance, while the LP relaxation yields bounds of poor quality. The VNS/TS heuristic obtains the optimal solution for $9$ out of $12$ instances with the largest gap being $4.30\%$.

As the number of customers increases to 10, the performance of CPLEX worsens significantly where only $1$ out of the $12$ instances are solved to optimality within the time and memory limits. {CPLEX is able to find feasible solutions for $11$ out of $12$ instances before termination yet the gaps are on average 34.31\%, which is significantly higher than the average gap of the VNS/TS heuristic 0.85\%.} 
The Lagrangian relaxation obtains bounds for 9 out of the 12 instances and the gap compared to the upper bound that is found by the heuristic is no more than $4.45\%$. We note that solving the Lagrangian subproblem is computationally challenging for the instances with $10$ customers and no solution is obtained within the time limits for 3 out of the 12 tested instances. {The lower bounds obtained by the LP relaxation are far away from the optimal values, which partially explains the low efficiency of CPLEX on this problem.}

For the instances with 15 customers, both CPLEX and the Lagrangian relaxation fail to solve the majority of the instances. CPLEX finds a feasible solution for only 1 out of the 11 instances. VNS/TS identifies solutions for all the instances within 317 seconds on average.

The results for the VNS/TS heuristic on medium instances are shown in Table \ref{BoundEvaluation_medium_instance}. The VNS/TS heuristic is able to solve all the instances with less than $35$ customers within $2$ hours (20 iterations), {allowing EV operators to implement the proposed approach for real-world tasks that do not require dynamic EV dispatching. Nevertheless, the proposed VNS/TS heuristic contains a large number of MIPs, which hinders its scalability. Future research efforts could possibly be made to re-formulate the embedded MIPs as quadratic knapsack problems and implement dynamic programming methods to solve them more efficiently.}

%======================================================================
\section{Case Study - Online Grocery Delivery in Kitchener-Waterloo} \label{chap:UseCase}
%======================================================================

In this section, the proposed EVRPTW-TP model along with the VNS/TS solution heuristic are evaluated on a case study of an EV fleet providing online grocery delivery services, which attracted great attention amid the COVID-19 pandemic in the Kitchener-Waterloo (KW) region in Ontario, Canada.
Through computational experiments, we particularly investigate the impact of electricity pricing schemes, time windows, winter range reductions, length of charging/discharging time period and fleet size on the fleet's routing and charging/discharging behaviors based on which we provide managerial insights.

%======================================================================
\subsection{Case Setup} \label{sec:BasicInfo}
%======================================================================
% 
% \subsubsection{Online Grocery Delivery}

We consider a local grocery store that uses EVs to provide online grocery delivery services. This business model has been increasingly popular in recent years \citep{McKinsey2020Grocery}, especially during the COVID-19 pandemic where people avoided going to supermarkets and local stores.

In the KW region, retail giants, such as Walmart, Zehrs, and T\&T Supermarket, as well as local small businesses are all providing such services.
Customers would browse available items online and place the orders. Delivery providers perform touch-free delivery during the predetermined time slot. 
In this case study, we only consider orders that are scheduled for delivery on the next day and accordingly apply EVRPTW-TP to optimize the routing and charging/discharging of the EV fleet for the planning horizon of $5am$ - $12am$.
% 
% \subsubsection{EV Fleet Specifications}

For the purpose of the case study, we first consider a small fleet consisting of $3$ EVs. The fleet size is then increased to 6 in order to investigate the impact on fleet operations (Section \ref{sec:FleetSize}). Service is provided to $30$ customers in the designated area, which reflects a real situation for a small business. 
The EVs in the fleet are homogeneous with identical range and cargo capacity. The EV that is considered has a range of $150$ kilometers with a $32.4$ kWh battery and a cargo capacity of $200$ order units. The travel time between locations is calculated based on Google Maps.
% 
% \subsubsection{Customers, Stations and the Depot}

\begin{figure}
    \centering
    \includegraphics[width=15cm,keepaspectratio]{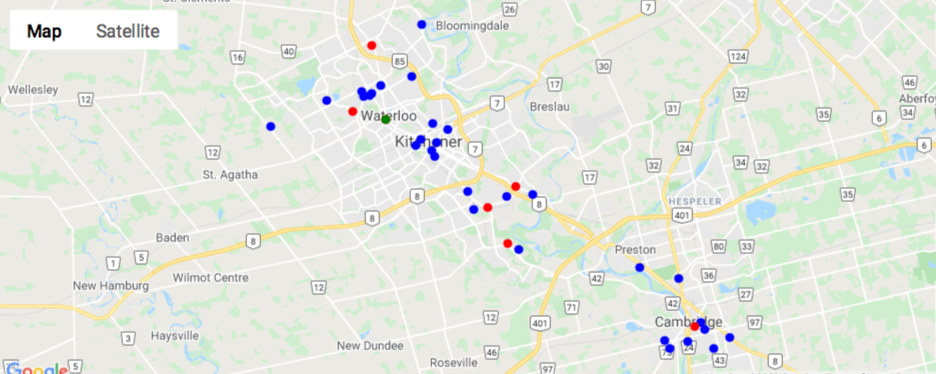}
    \caption{Selected Customers and Stations}
    \label{fig:vertex_location}
\end{figure}

The geographical locations of the actual customers, the stations, and the depot are shown in Figure \ref{fig:vertex_location} as blue, red, and green dots, respectively. The planning horizon is first divided into 3 time periods: morning ($5am$ - $12pm$), afternoon ($12pm$ - $6pm$), and evening ($6pm$ - $12am$). Results with different time windows are presented in Section \ref{sec:TimeWindow}. The service time at each customer location is set to $10$ minutes. The 7 selected stations are level-2 stations which are compatible with most EV models. The power supply of a level-2 charger is $240VAC/30A$. 
One hour of charge using a level-2 station adds $30$ kilometers of range \citep{Ministry2020} for the EV model that is considered in this case, i.e. $B = 270$ and $\alpha = 1.8$. 

{The solution time for the cases is on average 45 minutes which allows the fleet operator to set the routing, charging and discharging schedule at night. }

%======================================================================
\subsection{The Impact of the Electricity Pricing Scheme} \label{sec:Pricing}
%======================================================================
In this section, we study the impact of the electricity pricing on the EV fleet operations. As shown in Table \ref{tab:pricing_schemes}, we consider four electricity pricing schemes: 
\begin{description}
 \item[Scheme A:] The price varies between on-peak, mid-peak, and off-peak hours.
 \item[Scheme B:] The reward for energy discharging is the same as the price for energy charging.
 \item[Scheme C:] The price is only lower in the evening hours.
 \item[Scheme D:] Only energy charging is allowed (energy discharging to the grid is not allowed).
 \end{description}
 
\begin{figure}
    \includegraphics[width=15cm,keepaspectratio]{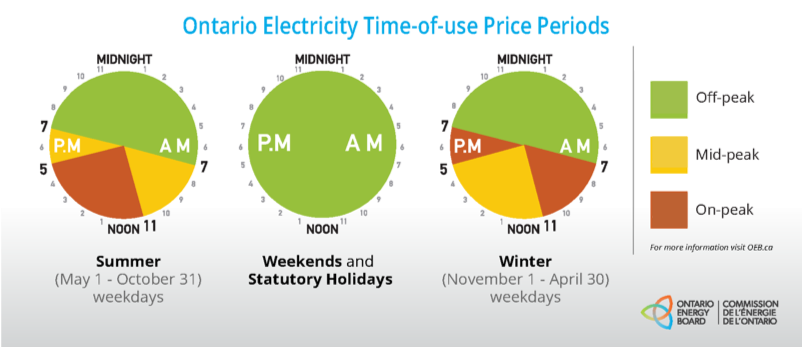}
    \caption{Time-of-Use Electricity Price Periods}
    \label{fig:tou_period}
\end{figure}

In scheme A, the charging rates are the real hydro rates that were in effect between May 1, 2018 and Nov 1, 2019. We note that the province of Ontario currently does not have a set reward for injecting electricity to the grid. Thus we evaluate the case where the reward for discharging is lower than the price of electricity (Scheme A) and the case where the reward is the same as the price of energy (Scheme B). We also evaluate the case where discharging is not allowed (Scheme D).
The electricity price periods are defined by the Ontario Energy Board differently for the summer and the winter (see Figure \ref{fig:tou_period}). We thus evaluate both summer and winter periods.  

\begin{table}
\centering
\caption{Electricity Prices}
\scriptsize
\begin{tabular}{@{}ccclcclcc@{}}
\toprule
       & \multicolumn{2}{c}{On-Peak ($\cent$/kWh)} &  & \multicolumn{2}{c}{Mid-Peak ($\cent$/kWh)} &  & \multicolumn{2}{c}{Off-Peak ($\cent$/kWh)} \\ \cmidrule(lr){2-3} \cmidrule(lr){5-6} \cmidrule(l){8-9} 
Scheme & Charge          & Discharge         &  & Charge          & Discharge          &  & Charge          & Discharge          \\ \midrule
A      & 13.40           & 10.00             &  & 9.40            & 8.00               &  & 6.50            & 6.50               \\
B      & 13.40           & 13.40             &  & 9.40            & 9.40               &  & 6.50            & 6.50               \\
C      & 9.40            & 8.00              &  & 9.40            & 8.00               &  & 6.50            & 6.50               \\
D      & 13.40           & 0.00              &  & 9.40            & 0.00               &  & 6.50            & 0.00               \\ \bottomrule
\end{tabular}
\label{tab:pricing_schemes}
\end{table}

The results obtained by solving EVRPTW-TP under the different electricity pricing schemes are shown in Table \ref{tab:PriceSchemeAnalysis}. 
The first two columns present the total distance traveled and the total electricity cost of the fleet. The total cost is decomposed into charging cost and discharging reward illustrated in the following two columns. The negative numbers in the ``Total Cost'' column indicate that the EV fleet makes a net profit due to selling energy. The last three columns present the number of hours spent by the fleet on charging and discharging electricity while on route and the number of station visits of the fleet, respectively.
{
\begin{table}
\centering
\caption{Overall Electricity Cost and Charging/Discharging Hours under Pricing Schemes}
\scriptsize
\begin{tabular}{@{}cccccccc@{}}
\toprule
Scheme     & Distance (km) & Total Cost ($\cent$)      & Char. Cost ($\cent$)     &  Disch. Reward ($\cent$)  & Char. (hr)    & Disch. (hr)      & \# of Station Visits \\ \midrule
A (summer) & 154.89        & -9.36                     & 638.64                   &  648.00                   & 2             & 9                & 3                    \\
A (winter) & 143.42        & -29.74                    & 690.26                   &  720.00                   & 3             & 10               & 3                 \\
B (summer) & 187.29        & -233.86                   & 730.94                   &  964.80                   & 3             & 10               & 3                   \\
B (winter) & 180.75        & -279.72                   & 878.04                   &  1157.76                  & 5             & 12               & 3                  \\
C          & 148.24        & 132.55                    & 535.75                   &  403.20                   & 0             & 7                & 0                  \\
D          & 143.88        & 202.03                    & 202.03                   &  0.00                     & 0             & 0                & 0                  \\ \bottomrule
\end{tabular}
\label{tab:PriceSchemeAnalysis}
\end{table}
}

The first insight is that allowing energy discharging is economically positive to EV fleet owners. 
When discharging is not allowed (Scheme D), the fleet's daily energy cost is $\$2.02$, while the gain reaches $\$2.80$ per day (Scheme B, winter) when discharging is allowed. The daily profits due to energy discharging could be as large as \$4.82 (the difference between schemes B and D), corresponding to an annual profits of $\$1,759.30$. According to \citet{fernandez2013capacity} and \citet{millner2010modeling} the EV battery life span is around $4.6$ - $5.5$ years under 100\% depth of discharge (DoD) daily usage in mild temperature. Therefore, considering the 8-year warranty offered by most EV producers, then enabling discharging can result in $2.5$ - $3.4$ years of reduction in battery life span. The total profits the fleet is able to gain from discharging before its end of life is \$7,978 - \$9,539.  Assuming EV battery price of 163 \$/kWh \citep{scott2020Ever}, the battery degradation cost is \$4,951 - \$6,734 for a fleet of $3$ EVs, each with a battery capacity of $32.4$ kWh. The profits gained from discharging can cover the battery degradation cost. Considering that the EV drivers are paid a fixed monthly salary, the additional profits assist to reduce the EV fleet's operational cost.
 %  considering the 8-year warranty offered by most EV producers, --> is this needed??
 
\begin{figure}[h]
    \centering
    \begin{minipage}[t]{0.48\textwidth}
        \includegraphics[width=3in]{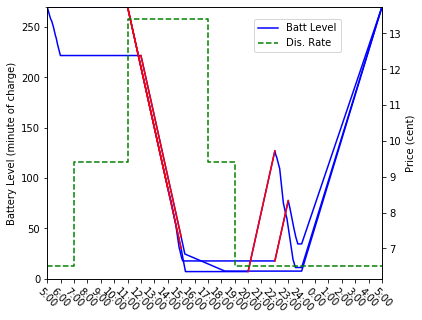}
        %\caption{Battery level (summer)}
        %\label{fig:batt_B_summer}
    \end{minipage}
    \begin{minipage}[t]{0.48\textwidth}
        \includegraphics[width=3in]{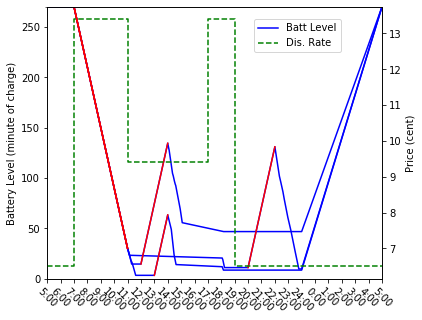}
        %\caption{Battery level (winter)}
        %\label{fig:batt_B_winter}
    \end{minipage}
    \caption{EVs' Battery Levels through a Day ($K = 3$, left: summer, right: winter)}
    \label{fig:BattLvlSchemeB}
\end{figure}

Not surprisingly, electricity prices impact fleet operations. 
In particular, since the reward rate in scheme B is much higher than in scheme A, EVs are more willing to discharge during peak hours and detour to recharge their batteries later. 
The travel distance under scheme B is, on average $32$ kilometers longer than under scheme A. 
EVs also spend $1$ to $3$ more hours discharging energy under scheme B than under scheme A. 
Moreover, the overall electricity cost varies significantly between these two cases. For schemes A and B, the total electricity cost is lower in winter than in summer. 
Taking scheme B as an example, in Figure \ref{fig:BattLvlSchemeB} each blue line depicts one EV battery level through the day and the charging/discharging behaviors are highlighted in red. The green dashed lines represent the discharging reward rates. In summer, in order to perform discharging during peak hours ($11am$ - $5pm$), at least one EV has to detour during the day because some customers should be served before $12pm$ (the end of the morning period). However, in winter, since the on-peak periods are $7am$ - $11am$ and $5pm$ - $7pm$ all the EVs can stay at the depot before $11am$ to perform discharge and leave to perform delivery services after that. As a consequence, the travel distance in winter is shorter than in summer for both schemes. Accordingly, EV operators can adjust their time slot settings with respect to the season. For instance, reducing the morning service period may decrease the overall electricity cost. Accounting for the changes in time-of-use prices and the corresponding time between the seasons is therefore important when designing the service time slots. Thus, next we evaluate the impact of the time windows on the EV operations.

%======================================================================
\subsection{The Impact of Time Windows} \label{sec:TimeWindow}
%======================================================================

In this section, we consider three variations of the time windows. We first consider 2-hour time windows, then the 3-periods presented in the previous section, and finally eliminate the time windows. Currently, 2-hour time windows are the most commonly used for grocery delivery in the KW area to allow customers to know the time that they would be served more precisely, thus reducing their waiting times and enhancing user experience. This of course comes at the expense of the flexibility of the EV fleet in terms of routing and scheduling. 

\begin{figure}
\centering
    \begin{minipage}[t]{0.9\textwidth}
    \centering
        \includegraphics[width=6in]{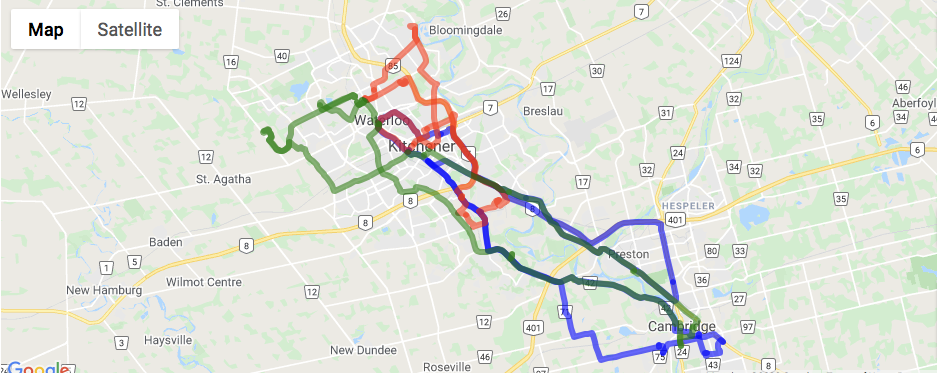} 
        %\caption{Routes (2-hour, winter)}
        %\label{fig:route_2hour_W}
    \end{minipage}
    
    \begin{minipage}[t]{0.9\textwidth}
    \centering
        \includegraphics[width=6in]{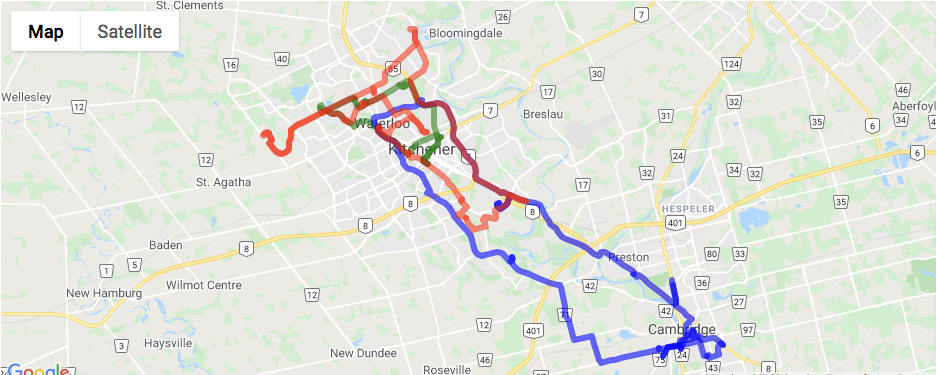} 
        %\caption{Routes (3-period, winter)}
        %\label{fig:route_3period_W}
    \end{minipage}
    
    \begin{minipage}[t]{0.9\textwidth}
    \centering
        \includegraphics[width=6in]{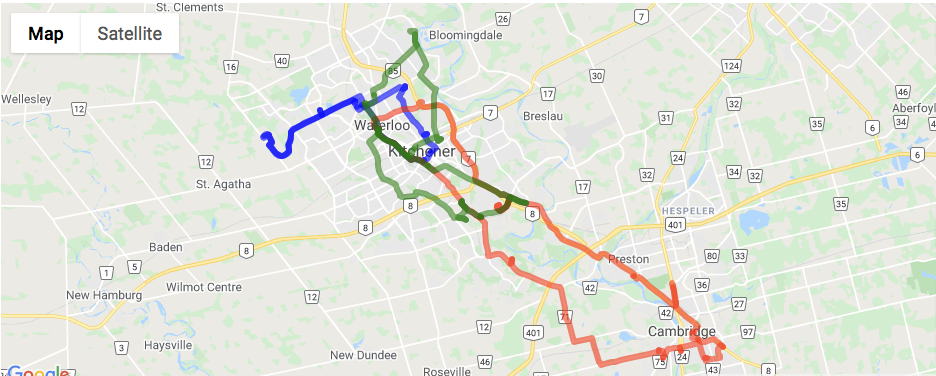} 
        %\caption{no-tw (winter)}
        %\label{fig:route_notw_W}
    \end{minipage}
    \caption{EV Routes (from top to bottom: 2-hours, 3-periods and no-tw)}
    \label{fig:routes_winter}
\end{figure}

The results presented in Table \ref{tab:timeWindowAnalysis} and Figure \ref{fig:routes_winter} reflect the trade-off between timely service and operational flexibility.
Under the 2-hours setting, customers that are geographically close might book time slots that are very far away from each other. 
In that case, an EV cannot always wait at a customer's location until the beginning of the time slot of the nearest customer. 
Instead, the EV will travel to a customer that is relatively distant but with a time slot starting earlier, leaving the nearest customer to be served by another EV. When a 3-period time window is used, more EVs are able to serve the customers in an order based on their geographical locations. As shown in Figure \ref{fig:routes_winter}, two EVs have to travel to Cambridge under the 2-periods setting, while only one has to do so under the other two settings. The increased flexibility leads to a decrease in the total travel distance from over $210$ kilometers in the 2-hours setting to around $150$ kilometers in the 3-period setting. The driving distance decreases to $131.81$ kilometers for the case of no time window constraints. 

In terms of charging and discharging activities, intuitively the fleet has the greatest flexibility in scheduling discharging in the case where no time windows are used. In both winter and summer, the fleet is able to discharge for $12$ hours in total during the planning horizon. For the ``3-period'' case, the fleet performs 9 hours of discharging, while for the 2-hour time windows, the fleet performs $6$ hours of discharging in summer, and $7$ hours in winter. 
As a consequence, as shown in Table \ref{tab:timeWindowAnalysis}, the cost of electricity is $1.85$ dollars and $1.64$ dollars for summer and winter respectively under the 2-hours setting. A profit of $1.17$ dollars and $1.12$ dollars is generated when no time windows are enforced in summer and winter, respectively. This difference can be seen as the price of providing ``timely service'' delivery.
{
\begin{table}
\centering
\caption{Overall Electricity Cost under Time Window Settings}
\scriptsize
\label{tab:timeWindowAnalysis}
\begin{tabular}{@{}ccccccccc@{}}
\toprule
Season                  & TW Setting & Dist. (km) & Total Cost ($\cent$) & Char. Cost ($\cent$)     &  Disch. Reward ($\cent$)   & Char. (hr) & Disch. (hr) & \# of Stat. Visits \\ \midrule
\multirow{3}{*}{Summer} & 2-hours     & 229.19        & 184.90              & 602.50                   & 417.60                     & 0             & 6             & 3  \\
                        & 3-periods   & 154.89        & -9.36               & 638.64                   & 648.00                     & 2             & 9             & 3  \\
                        & no-tw       & 131.81        & -117.36             & 746.64                   & 864.00                     & 4             & 12            & 3  \\ \midrule
\multirow{3}{*}{Winter} & 2-hours     & 212.78        & 164.09              & 596.09                   & 432.00                     & 2             & 7             & 3  \\
                        & 3-periods   & 143.42        & -29.74              & 690.26                   & 720.00                     & 3             & 10            & 3  \\
                        & no-tw       & 135.95        & -111.53             & 752.47                   & 864.00                     & 4             & 12            & 3 \\ \bottomrule
\end{tabular}
\end{table}
}

%======================================================================
\subsection{The Impact of Range Reduction} \label{sec:ShorterRange}
%======================================================================

{
In this section, we investigate the impact of range reduction under low ambient temperature. According to the data presented by \citet{yuksel2015effects} and \citet{lohse2013ambient}, an EV's driving range could drop by 50\% in severe weather conditions. We thus consider EV ranges from 100\% to 50\% of its stated range as shown in Table \ref{tab:RangeAnalysis}. We note that battery capacity is assumed to be fixed, we vary EV range by adjusting the energy consumption rate. For example, increasing energy consumption rate by $100\%$ is equivalent as reducing EV range by $50\%$. We also assume the charging and discharging processes are not affected by ambient temperature.}

{
Not surprisingly, Table \ref{tab:RangeAnalysis} shows that the fleet's total energy cost increases as EV range drops. When EVs are operated at close-to-full ranges ($90\%$ or more), the fleet can discharge for $10$ hours during peak periods and still have the capacity to visit charging stations 3 times to charge the energy back during off-peak hours. Consequently, the discharging rewards cancel out all the charging costs. When EV range drops to $70$-$80\%$, an EV can not make its detour for en-route charging to have the same amount of energy being discharged during peak periods. It thus discharges one hour less and cancels the en-route station visit, which leads to higher total costs but less driving distance. If the EV range drops to $50$-$60$\% of its stated range, the fleet is still able to discharge for $8$ hours, yet the total cost increases to around $2$ dollars per day as more charging is needed.}

\begin{table}
\centering
{
\caption{Overall Electricity Cost and Charging/Discharging Hours under Range Reductions in Winter}
\scriptsize
\begin{tabular}{@{}cccccccc@{}}
\toprule
Range & Distance (km) & Total Cost ($\cent$) & Char. Cost ($\cent$) & Disch. Reward ($\cent$) & Char. (hr) & Disch. (hr) & \# of Station Visits \\ \midrule
100\%            & 143.42        & -29.74               & 690.26               & 720.00                     & 3          & 10          & 3 \\
90\%             & 145.81        & -3.65                & 716.35               & 720.00                     & 3          & 10          & 3 \\
80\%             & 141.38        & 21.32                & 669.32               & 648.00                     & 3          & 9           & 2 \\
70\%             & 141.38        & 56.76                & 704.76               & 648.00                     & 3          & 9           & 2 \\
60\%             & 142.50     & 131.83               & 707.83               & 576.00                     & 3          & 8           & 3 \\
50\%             & 142.58        & 198.78               & 774.78               & 576.00                     & 4          & 8           & 3 \\ \bottomrule
\end{tabular}
}
\label{tab:RangeAnalysis}
\end{table}

%======================================================================
\subsection{The Impact of Time Period Length} \label{sec:TimePeriod}
%======================================================================
In this section, we investigate the impact of the length of each charging/discharging time period. Table \ref{tab:PeriodLengthAnalysis} presents the fleet's energy cost, travelling distance, and the number of station visits when each time period lasts $15$, $30$, and $60$ minutes under winter and summer hours, respectively.

\begin{table}[H]
\centering
\caption{Overall Electricity Cost under Different Time Period Lengths}
\scriptsize
\label{tab:PeriodLengthAnalysis}
\begin{tabular}{@{}ccccccccc@{}}
\toprule
Season                  & Perid Length ($\delta$) & Dist. (km) & Total Cost ($\cent$) & Char. Cost ($\cent$)     &  Disch. Reward ($\cent$)   & Char. (hr) & Disch. (hr) & \# of Stat. Visits \\ \midrule
\multirow{3}{*}{Summer} & 60     & 154.89        & -9.36               & 638.64                   & 648.00                     & 2             & 9             & 3  \\
                        & 30     & 145.41        & -47.84              & 672.16                   & 720.00                     & 1             & 10            & 3  \\
                        & 15     & 145.19        & -48.15              & 707.85                   & 756.00.                    & 1.5           & 10.5          & 3 \\ \midrule
\multirow{3}{*}{Winter} & 60     & 143.42        & -29.74              & 690.26                   & 720.00                     & 3             & 10            & 3  \\
                        & 30     & 148.67        & -60.19              & 767.81                   & 828.00                     & 3             & 11.5          & 3  \\
                        & 15     & 149.14        & -65.83              & 780.17                   & 846.00                     & 3             & 11.75         & 3 \\ \bottomrule
\end{tabular}
\end{table}

Shortening time periods enhance the fleet's operational flexibility, allowing fleet operators to consider schedules of higher resolution. For instance, consider an EV that visits a station at 11:15am and plan to discharge for one hour during the visit, the EV would have to wait for $45$ minutes before it starts discharging if $\delta=60$, while the waiting time is $15$ and $0$ minutes when $\delta=30$ and $\delta=15$, respectively. Due to the reduced waiting time, the EVs can better coordinate their charging/discharging and routing schedules and spare more time for grid ancillary services. Moreover, we observe that some EVs may finish serving assigned customers and return to the depot before the end of evening peak hours. These EVs usually don't have enough energy to discharge for one hour, thus can not make additional profits when $\delta=60$. Nevertheless, when $\delta=15$, the fleet is allowed to inject less energy back to the grid each time, hence can jointly discharge for $1.5$ and $1.75$ additional hours under summer and winter settings, respectively. Consequently, as presented in Table \ref{tab:PeriodLengthAnalysis}, the total energy cost decreases as the value of $\delta$ decreases.

Although shorter time periods lead to decreased energy cost, it also introduces extra computational burdens. The average computational time for the cases of $\delta=60$, $30$ and $15$ is $45$, $110$, and $160$ minutes, respectively. In practice, fleet operators can choose the length of time period based on the available computational resources and the problem size.

%======================================================================
\subsection{The Impact of EV Fleet Size} \label{sec:FleetSize}
%======================================================================

\begin{figure}[th]
\centering
    \centering
        \begin{minipage}[t]{0.48\textwidth}
        \includegraphics[width=3in]{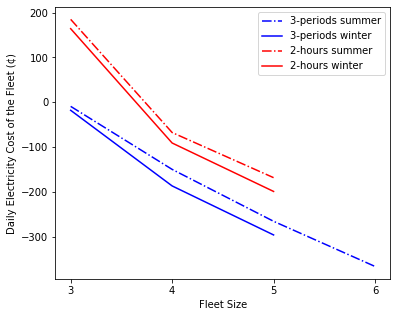} 
        \caption{Electricity Cost under Different Fleet Sizes}
        \label{fig:fleet_size_analysis}
    \end{minipage}
    \begin{minipage}[t]{0.48\textwidth}
        \includegraphics[width = 3in]{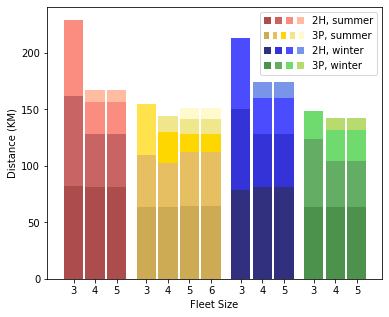}
        \caption{Overall Travelling Distance under Different Fleet Sizes}
        \label{fig:fleet_size_analysis_distance}
    \end{minipage}
\end{figure}

In this section, we consider the impact of EV fleet size on the overall electricity cost under the two restrictive time window settings, 2-hours and 3-periods. The number of EVs is varied from 3 to 6 and the corresponding results are shown in Table \ref{tab:FleetSizeAnalysis} and Figure \ref{fig:fleet_size_analysis}.
Intuitively, increasing fleet size reduces the overall travel distance and provides more spare time to discharge energy to the grid.

\begin{figure}
\centering
    \begin{minipage}[t]{0.48\textwidth}
    \centering
        \includegraphics[width=3in]{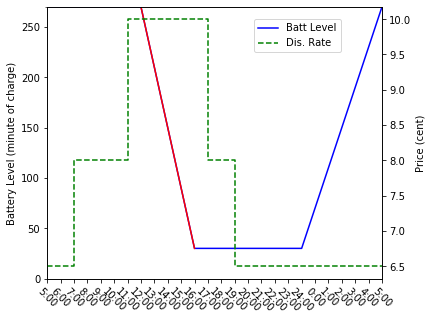} 
        %\caption{Routes (2-hour, winter)}
        %\label{fig:soley_summer}
    \end{minipage}
    \begin{minipage}[t]{0.48\textwidth}
    \centering
        \includegraphics[width=3in]{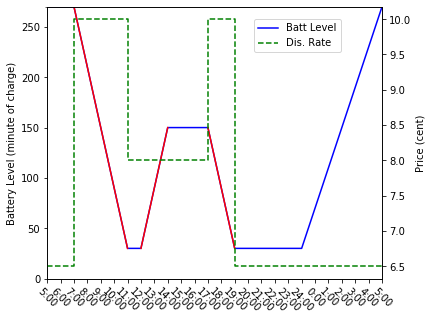} 
        %\caption{Routes (3-period, winter)}
        %\label{fig:soley_winter}
    \end{minipage}
    \caption{Battery Level of an EV without Logistic Service (left: summer, right: winter)}
    \label{fig:soley_ancillary}
\end{figure}

However, given the high acquisition cost, using an EV solely for the purpose of trading energy with the grid is currently not economical. Figure \ref{fig:soley_ancillary} shows the battery levels of an EV during summer and winter hours. In summer, the EV discharges its battery during the on-peak rate for 4 hours and recharges under the off-peak rate at night, bringing a profit of $1.00$ dollar. 
Similarly in winter, the EV makes a profit of $1.09$ dollars by discharging for 6 hours under on-peak rates and recharging for 2 and 4 hours under mid-peak and off-peak rates respectively. Adding one additional EV increases the profit by $1.01$ dollars in summer hours an $1.09$ dollars in winter hours. In that case the additional EV is acting as a static battery storage.

\begin{table}
\centering
\caption{Overall Electricity Cost of EV fleets of Different Sizes}
\label{tab:FleetSizeAnalysis}
\begin{tabular}{@{}cccccc@{}}
\toprule
\multicolumn{1}{l}{}       & \multicolumn{1}{l}{} & \multicolumn{4}{c}{Number of EVs} \\ \cmidrule(l){3-6} 
TW                         & Season               & 3      & 4      & 5      & 6      \\ \midrule
\multirow{2}{*}{2-hours}   & Summer               & 184.90  & -67.75 & -168.55 & -      \\
                           & Winter               & 164.09  & -91.15 & -199.30 & -      \\  \midrule
\multirow{2}{*}{3-periods} & Summer               & -9.36  & -150.05 & -265.97 & -366.77 \\
                           & Winter               & -29.74  & -186.98 & -296.42 & - \\
\bottomrule
\end{tabular}
\end{table}

The impact of fleet size on the travel distance is shown in Figure \ref{fig:fleet_size_analysis_distance}. Each bar consists of several sub-bars highlighted in different colors each representing the travel distance of an EV in the fleet. As shown, for all the four cases, when the fleet size is $3$, the fleets' overall travel distance is the highest. As the fleet size increases, the overall travel distance decreases. For the ``2-hours'' time window cases (red and blue bars), increasing the fleet size from $3$ to $4$ results in a relatively large reduction in the overall travel distance, while for other cases, the distance reductions are marginal. This observation suggests that the newly added EVs contribute to reducing the overall electricity cost mainly through grid ancillary services instead of improving the routing. Additional EVs thus are increasingly playing the role of a battery storage which is not economical under current energy prices.

\section{Conclusion}  \label{chap:Conclusion}

This paper presented the EV routing problem with time windows under time-variant electricity prices. The formulated problem jointly optimizes the routing of a fleet of electric vehicles along with the scheduling of the charging and discharging activities. The proposed framework permits the EV fleet to provide the primary service of delivering products to customers while minimizing the energy costs by scheduling energy charging while taking into consideration the varying energy prices throughout the planning horizon. Furthermore, the electric vehicles can store energy at periods with relatively low energy prices and inject it back to the grid at high price periods which helps in cost recovery and potentially realize profits. This framework demonstrates the capabilities of EV fleets in helping energy grids to smooth out the demand and reduce the gap between on-peak and off-peak by acting as energy storage. The proposed problem is formulated as a multiperiod optimization model that is challenging to solve to optimality. A Lagrangian relaxation that provides lower bounds is proposed. To find feasible solutions quickly which is important for implementation in practice, a hybrid variable neighborhood search/tabu search heuristic which finds high quality solutions in relatively short amount of computational time is presented. Numerical results on commonly used instances form the literature demonstrated the capabilities of the proposed heuristic in finding high quality solutions. A case study on a grocery delivery service from the region of Kitchener-Waterloo in Ontario Canada further illustrates the impact of energy prices, time windows design, winter range reduction, time period length and the size of the EV fleet on operations and costs.

The proposed framework is a step towards realizing the full integration of transportation and energy networks. The joint optimization of both systems promises increasing opportunities to improve reliability and reduce costs. {Nevertheless, it is not without limitations. Some assumptions, such as linear charging/discharging, constant energy consumption, are made to have a tractable problem. Moreover, battery degradation cost which could potentially be incurred by frequent charging and discharging is not taken into account as it leads to longer time horizon than the one considered. Incorporating these realistic characteristics presents interesting future research directions. From the perspective of solution method, the proposed Lagrangian relaxation approach and the VNS/TS heuristic both have sub-problems that are computationally demanding, which hinders their scalability. Future research efforts could be made to improve solution efficiency of both approaches.}

While in the framework that is proposed in this paper, the energy prices in a period are known in advance, a model that integrates dynamic and uncertain prices of energy will be investigated in future work. Furthermore, future work will particularly focus on integrating energy generated from renewable sources to further advance the development of green logistics systems.
\section*{Acknowledgment}

Bo Lin was supported by the Energy Council of Canada energy policy research fellowship and Bissan Ghaddar was  supported by NSERC Discovery Grant 2017-04185. We greatly thank the referees, the Associate Editor, and the Editor for their thorough and thoughtful comments that helped us improve the quality of the paper.

% Acknowledgments here
% \ACKNOWLEDGMENT{%
% Enter the text of acknowledgments here
% }% Leave this (end of acknowledgment)

% Appendix here
% Options are (1) APPENDIX (with or without general title) or 
%             (2) APPENDICES (if it has more than one unrelated sections)
% Outcomment the appropriate case if necessary
%
% \begin{APPENDIX}{<Title of the Appendix>}
% \end{APPENDIX}
%
%   or 
%
% \begin{APPENDICES}
% \section{<Title of Section A>}
% \section{<Title of Section B>}
% etc
% \end{APPENDICES}

% References here (outcomment the appropriate case) 

% CASE 1: BiBTeX used to constantly update the references 
%   (while the paper is being written).
%\bibliographystyle{ormsv080} % outcomment this and next line in Case 1
%\bibliography{<your bib file(s)>} % if more than one, comma separated

% CASE 2: BiBTeX used to generate mypaper.bbl (to be further fine tuned)
%\input{mypaper.bbl} % outcomment this line in Case 2

% \bibliographystyle{ormsv080}
\bibliographystyle{abbrvnat}
\bibliography{paper.bbl}

\end{document}